\documentclass[12pt]{amsart}

\usepackage{amsfonts,amsmath} 
 
\usepackage{amssymb,latexsym}

\newtheorem{theorem}{Theorem}[section]  
\newtheorem{lemma}[theorem]{Lemma}  
\newtheorem{proposition}[theorem]{Proposition}

\newtheorem*{claim}{Claim}  
\newtheorem*{conjecture}{Conjecture}

\newtheorem{definition}{Definition}[section]

\theoremstyle{remark}  
\newtheorem{remark}{Remark}[section]

\numberwithin{equation}{section}

\newcommand{\beq}{\begin{equation}}  
\newcommand{\enq}{\end{equation}}  
  
\newcommand{\ltwo}{L^{2}}  
\newcommand{\Sc}{\mathcal{S}}

\newcommand{\rn}{\mathbb{R}^{n}}

\newcommand{\real}{\mathbb{R}}  
\newcommand{\R}{\mathbb{R}}

\newcommand{\N}{\mathbb{N}}

\newcommand{\oOmega}{\vartheta}  
\newcommand{\loc}{\mbox{\footnotesize{loc}}}  
\newcommand{\cal}{\mathcal}   
\renewcommand{\Hat}{\widehat} 

\DeclareMathOperator{\supp}{Supp}  
\DeclareMathOperator{\dive}{div}

\begin{document}  
\title[Weak solutions and enstrophy defects]{Weak solutions,  
renormalized solutions and enstrophy defects in 2D turbulence}  
  
\author[Milton Lopes {\em et al.}]{Milton C. Lopes Filho}  
\address{Depto. de Matem\'{a}tica, IMECC-UNICAMP \\ Cx. Postal 6065, Campinas SP 13083-970, Brazil} 
\email{mlopes@ime.unicamp.br}  
  
\author[]{Anna L. Mazzucato}  
\address{Department of Mathematics, Penn State University \\ 218
  McAllister Building, University Park, PA 16802, U.S.A.}  
\email{mazzucat@math.psu.edu}  
  
\author[]{Helena J. Nussenzveig Lopes}  
\address{Depto. de Matem\'{a}tica, IMECC-UNICAMP \\ Cx. Postal 6065, Campinas SP 13083-970, Brazil}
\email{hlopes@ime.unicamp.br}

\begin{abstract} 

 Enstrophy, half the integral of the square of vorticity, plays a role in 2D turbulence theory analogous to
that played by kinetic energy in the Kolmogorov theory of 3D
turbulence. It is therefore interesting to
obtain a  description of the way enstrophy is dissipated at high Reynolds number. In this article we explore the 
notions of viscous and transport enstrophy defect, which model the spatial structure of the dissipation of enstrophy.
These notions were introduced by G. Eyink in an attempt to reconcile the Kraichnan-Batchelor theory of 2D 
turbulence with current knowledge of the properties of weak solutions of the equations of incompressible and 
ideal fluid motion. Three natural questions arise from Eyink's theory: (1) Existence of the enstrophy 
defects (2) Conditions for the equality of transport and viscous enstrophy defects (3) Conditions for 
the vanishing of the enstrophy defects. In \cite{Ey}, Eyink proved a number of results related to these 
questions and formulated a conjecture on how to  answer these problems in a physically meaningful context. 
In the present article we improve and extend some of Eyink's results and present a counterexample to his conjecture. 
\end{abstract} 
 
\maketitle   
{\it \small Mathematics Subject Classification: 76F02, 35Q35, 76B03} 
 
{\it \small Keywords: Incompressible flow, two-dimensional flow, turbulence, vorticity} 
\tableofcontents 
  
\section{Introduction}  
 
This article is concerned with certain properties of weak solutions of the incompressible 
Euler equations in two space dimensions and with the corresponding vanishing viscosity limit
in connection with the modeling of two-dimensional turbulence. To
put our discussion in context it is useful to recall some of the basic features of the
Kraichnan-Batchelor (KB)  theory of two-dimensional turbulence, introduced in 
\cite{Kraich,Ba}. This is a phenomenological theory,  modeled 
after Kolmogorov's theory of 3D turbulence. The notion of {\em enstrophy} cascade plays a 
central role in KB theory, similar to the role of the energy cascade in Kolmogorov's 
theory. Enstrophy is half the integral of the square of vorticity, a conserved 
quantity for smooth ideal 2D flow, which is dissipated in viscous flow. In the cascade 
picture, the nonlinearity transports enstrophy from large to small scales, 
where it is dissipated by viscosity. A key issue in the KB theory is that
such a picture must be sustained as viscosity vanishes, in a way that allows 
the rate at which enstrophy is dissipated to remain bounded away from zero as viscosity 
disappears. For details and the associated literature we refer the reader to \cite{Frisch}, 
especially Section 9.7, and references there contained. 

Let us consider a family of viscous flows, which we assume to have uniformly bounded enstrophy
as viscosity vanishes. This sequence is compatible with the KB cascade if the enstrophy
dissipation rate is bounded away from zero. Taking subsequences as needed, such a  family leads to  a weak 
solution of the 2D incompressible  Euler equations, see \cite{Maj}, which must dissipate enstrophy. 
The difficulty one faces is that weak solutions of  the incompressible 2D 
Euler equations with finite enstrophy conserve enstrophy exactly, a known fact which 
we will examine in detail later. We note that this difficulty does not occur in 
3D, as energy dissipative solutions of the incompressible 3D Euler equations with finite initial energy 
have been shown to exist, see \cite{DR,Shni}. 

        Recently, G. Eyink proposed a way around the paradox     
outlined above, see \cite{Ey}, by considering flows with unbounded 
local enstrophy. Eyink's idea raises the mathematical problem  of assigning meaning to enstrophy 
dissipation for flows with infinite enstrophy. In \cite{Ey}, Eyink introduced two notions of 
enstrophy defect in his attempt to describe the spatial structure of the enstrophy dissipation. 
These enstrophy defects are limits of enstrophy source terms in approximating enstrophy balance equations. 
When the relevant approximation is vanishing viscosity, this limit gives rise to a viscous enstrophy defect. 
The other defect introduced by Eyink was a purely inviscid enstrophy defect associated with mollifying a 
weak solution, which we call transport enstrophy defect. Eyink formulated a conjecture stating that both 
enstrophy defects are well-defined, that they give rise to the same
distribution in the limit and that they do not always 
vanish. One of the main purposes of the present work is to present a counterexample to Eyink's conjecture. 

       Beyond the description of 2D turbulence, there are two other concerns that motivate this paper. 
The first is the problem of uniqueness of weak solutions for incompressible 2D Euler, a 
long-standing open problem. Existence of weak solutions is known for compactly supported initial vorticities 
in the space $({\cal BM}_+ + L^1) \cap H^{-1}_{\loc}$, where ${\cal BM}_+$ is the cone of nonnegative bounded
Radon measures, see \cite{DeL,Sch,VW93}. In contrast, uniqueness of weak solutions is only 
known for vorticities which are bounded or nearly so, see \cite{Vi2,Yu63,Yu}. It is conceivable  that
the usual notion of weak solution is too weak to guarantee uniqueness, and that a criterion 
is required to select the `correct' weak solution. Properties that distinguish those weak solutions which are 
inviscid limits are particularly interesting, and we will encounter
some of these properties in this paper.  

The second concern is connected with the general issue of inviscid dissipation. Transport by smooth 
volume-preserving flows merely rearranges the transported quantity. This property is 
maintained even when the flow is not smooth, as long as we restrict ourselves to  renormalized solutions of 
the transport equations, in the sense of DiPerna and Lions, see \cite{DiPL}. 
Weak solutions (in the sense of distributions) of transport equations by divergence-free vector fields are
always renormalized solutions if the transported quantity and the transporting velocity are
sufficiently smooth. In the special case of weak solutions of the 2D Euler equations vorticity is always 
a renormalized solution of the vorticity equation, regarded as a linear transport equation,  as long as 
enstrophy is finite. Consequently, for finite enstrophy flows the distribution function of
vorticity is conserved in time. What happens with the distribution function of vorticity under less 
regular flows is a very interesting problem, closely related to the present work.    

           The remainder of this article is divided into six sections. In Section 2 we review the DiPerna-Lions  
transport theory and we apply it to ideal, incompressible, two-dimensional flow. In Section 3 we introduce the  
enstrophy defects, we prove that the viscous enstrophy defect vanishes for flows with finite enstrophy  
and we formulate a version of Eyink's conjecture. In Section 4 we  prove that the enstrophy density 
associated to a viscosity solution is a weak solution of a transport equation as long as vorticity lies 
in the space $L^2(\log L)^{1/4}$, an Orlicz space slightly smaller than $L^2$. We also show that the 
transport enstrophy defect exists as a distribution for vorticities in $L^2(\log L)^{1/4}$ and vanishes 
if the weak solution in this space happens to be an inviscid limit. 
In Section 5 we present examples showing that the results obtained in the previous section are nearly sharp.  
In Section 6 we exhibit a counterexample to Eyink's conjecture.  Finally, we draw some conclusions and highlight open
problems in Section 7. 

Technically speaking, we make use of the framework usually found in the study of nonlinear problems 
through weak convergence methods as well as  harmonic analysis and function space theory. One  distinction 
between our work and \cite{Ey} is that we consider flows in the full plane with compactly supported 
initial vorticity, whereas Eyink dealt with periodic flows. Working in the plane is convenient because of 
the simpler expression for the Biot-Savart law and because it is easier to find the function space results we 
require. The trade-off is the need  to work around problems arising from infinity, such as loss of tightness 
along  vorticity sequences.     
 
We conclude this introduction by fixing notation. We denote by $B(x;r)$ the disk centered at $x$ with radius $r$
in the plane. The characteristic function of a set $E$ is denoted by $\chi_E$. If $X$ is a function space then
$X_c$ denotes the subspace of functions in $X$ with compact support and $X_{\loc}$ denotes the space of functions
which are locally in $X$. We use alternatively $C^{\infty}_c$ or $\mathcal{D}$ to denote the space of smooth 
compactly-supported test functions. We use $W^{k,p}$ and $H^{s}$ to denote the classical Sobolev spaces. 
We denote by $L^{p,q}$ the Lorentz spaces and $B^s_{p,q}$ the Besov
spaces as defined respectively in \cite{BS88} and  \cite{B-L}.  
  
\section{Weak solutions and renormalized solutions}  
  
The purpose of this section is to discuss the relation between weak solutions of the incompressible   
2D Euler equations and DiPerna-Lions renormalized solutions of linear transport equations.   
  
We begin by recalling the vorticity formulation of the two dimensional Euler equations:  
\begin{subequations} \label{e:euler}  
 \begin{align}     
         &\partial_{t} \omega + u\cdot \nabla \omega=0, \label{e:euler.a}\\  
         &u=K\ast \omega, \label{e:euler.b}  
 \end{align}  
\end{subequations}  
with the Biot-Savart kernel $K$ given by 
\[ K(x) \equiv \frac{x^{\perp}}{2\pi|x|^2}, \]  
$(x_1,x_2)^{\perp} = (-x_2,x_1)$, and where the convolution in \eqref{e:euler.b}   
occurs only in the spatial variable.  Note that the specific form of the Biot-Savart kernel implies that 
$\dive u=0$.  
  
Identity \eqref{e:euler.a} is a transport equation for the vorticity.   
Therefore, if $u$ is sufficiently smooth so that $\omega$ is a classical   
solution, the vorticity itself and any function of it are transported along   
the flow induced by $u$. In particular, the enstrophy density function   
$\oOmega(x,t)=|\omega(x,t)|^{2}/2$ is conserved   
along particle trajectories, and, as the velocity $u$ is divergence-free,  
the enstrophy $\Omega (t)\equiv \int \oOmega (x,t)\, dx$   
is a globally conserved quantity in time.  
  
There is a well-developed theory of weak solutions for \eqref{e:euler}.  Well-posedness   
for weak solutions has been established for those initial vorticities which are bounded    
or nearly so, see \cite{Yu63,Yu,Vi1,Vi2}. If vorticity belongs to   
$L^p$ then, by Calderon-Zygmund theory and the   
Hardy-Littlewood-Sobolev inequality, $u \in W^{1,p}_{\loc}$ so that, if $p \geq 4/3$ then $u \in L^{p'}$ with   
$p'=p/(p-1)$. Hence the relevant nonlinear term, $u \, \omega$, is locally integrable and   
the transport equation \eqref{e:euler.a} lends itself to a standard weak formulation. To be precise we   
recall the weak formulation of the initial-value problem for \eqref{e:euler}. Let  
$\omega_0 \in L^p(\real^2)$, $p \geq 4/3$.   
  
\begin{definition} \label{weakvort}  
Let $\omega = \omega(x,t) \in L^{\infty}([0,T);L^p(\real^2))$ for some $p \geq 4/3$ and let $u = K \ast \omega$.   
We say $\omega$ is a weak solution of the initial-value problem for   
\eqref{e:euler} if, for any   
test function $\varphi \in C^{\infty}_c([0,T)\times \real^2)$, we have:  
\[ \int_0^T\int_{\real^2} \varphi_t \omega + \nabla \varphi \cdot u \, \omega \;dxdt +   
\int_{\real^2} \varphi(x,0)\omega_0(x)\;dx = 0.\]  
In addition, we require that the velocity field $u \in L^{\infty}([0,T);L^2(\real^2) + L^{\infty}(\real^2))$.  
\end{definition}  
  
Existence of weak solutions has been established for initial vorticities $\omega_0 \in ({\mathcal BM}_{c,+} + L^1_c)  
\cap H^{-1}_{\loc}$, see \cite{DiPM1,DeL,VW93,Sch}; however, these results require a more elaborate weak   
formulation in order to accommodate the additional irregularity in vorticity. If the vorticity is in $L^p$ for some  
$p \geq 4/3$ then all weak formulations reduce to the one in Definition \ref{weakvort}. In this paper we are mostly 
concerned with flows whose vorticity is in $L^2$ or nearly so, and for these flows,  Definition \ref{weakvort} is  
adequate. There is one situation of present interest for which Definition \ref{weakvort} cannot be used, namely, 
that of vorticities in the Besov space $B^{0}_{2,\infty}$. In this case a weak velocity formulation, see \cite{DiPM1}, 
should be used instead. 
  
Given that, for vorticities in $L^p$, the velocities are only $W^{1,p}_{\loc}$, it is natural to consider   
weak solutions of \eqref{e:euler} in the context of the  
theory of renormalized solutions for linear transport equations, introduced by  
DiPerna and Lions \cite{DiPL}. We recall below the definition of renormalized solution for linear transport   
equations without lower-order term.  
  
If $E \subseteq \rn$ then $|E|$ denotes the Lebesgue measure of $E$.  
Let $L^{0}$ be the set of all measurable functions $f$  
on $\rn$ such that $|\{|f(x)|>\alpha\}|<\infty$, for each $\alpha >0$.  
Let $v\in L^{1}([0,T];W^{1,1}_{\text{loc}})$ such that  
\begin{equation} \label{e:renorm0}  
      (1+|x|)^{-1}\,v \in L^{1}([0,T];L^{1}) + L^{1}([0,T];L^{\infty}).  
\end{equation}

\begin{definition}  
 A function $\omega \in L^{\infty}([0,T];L^{0})$ is called a {\em  
 renormalized} solution to the linear transport equation   
\[ \omega_t + v \cdot \nabla \omega = 0 \]  
if, in the sense of distributions,
\begin{equation} \label{e:renorm1}  
  \partial_{t} \beta(\omega) + v\cdot \nabla\beta(\omega) = 0,   
\end{equation}  
for all $\beta \in   
\mathcal{A}=\{\beta \in C^{1}, \beta \mbox{ bounded, vanishing near } 0 \} $.  
\end{definition}  
  
The most important property of renormalized solutions is that, in general, they are unique. The   
connection between weak solutions of the Euler equations and renormalized solutions of the vorticity  
equation \eqref{e:euler.a}, regarded as a linear transport equation with given velocity, is known.   
However, this relation has not been clearly stated in the literature. We address this omission in the following result.   
  
\begin{proposition} \label{weakrenorm}  
Let $p\geq 2$. If $\omega=\omega(x,t) \in L^{\infty}([0,T);L^p(\real^2))$ is a weak solution of the   
Euler equations then $\omega$ is a renormalized solution of transport equation \eqref{e:euler.a}   
with velocity $u = K\ast\omega$. Let $1<p<2$. If $\omega$ is a weak solution of the Euler equations obtained as   
a weak limit of a sequence of exact smooth solutions (generated, for example, by mollifying   
initial data and exactly solving the equations) then $\omega$ is a renormalized solution of \eqref{e:euler.a}.  
\end{proposition}  
  
\begin{proof}  
If $p\geq 2$, then the velocity $u$ belongs to $L^{\infty}([0,T);W^{1,p}_{\loc})$ and hence to   
$L^{\infty}([0,T);W^{1,p^{\prime}}_{\loc})$, as $p \geq p^{\prime}$.   
The velocity $u$ satisfies the mild growth condition \eqref{e:renorm0} because $L^2+L^{\infty}$ is  
contained in $L^1+L^{\infty}$ and an $L^2+L^{\infty}$ estimate on velocity was required in the definition 
of weak solution. Hence, we are under the conditions   
of the consistency result, Theorem II.3 in \cite{DiPL}, so we may conclude that $\omega$ is a   
renormalized solution. The statement regarding weak solutions that are limits of exact smooth solutions is a   
consequence of the stability result contained in Theorem II.4 in \cite{DiPL}.  
\end{proof}  
  
It is an interesting question whether the vanishing viscosity limit gives rise to a  
renormalized solution as well, if the initial vorticity is in $L^p$, $1<p<2$.   
  
Let $\omega \in L^{\infty}([0,T);L^2_c(\real^2))$ be a weak solution of \eqref{e:euler}. By Proposition \ref{weakrenorm}  
$\omega$ is also a renormalized solution. Since the velocity is divergence-free, we may conclude,   
using the full strength of the DiPerna-Lions theory of renormalized solutions, that the distribution function of $\omega$  
is time-independent, i.e.:  
\begin{equation} \label{distindep}  
\lambda_{\omega}(s,t) \equiv |\{x \in\real^2 \, | \, |\omega(x,t)|>s \}| =\lambda_{\omega}(s,0) \equiv   
\lambda_{\omega_0}(s),  
\end{equation}   
see the second Theorem III.2 of \cite{DiPL}. Therefore, all rearrangement-invariant norms of vorticity are conserved  
in time. In particular, the enstrophy $\Omega(t)$ is preserved for any weak solution of the 2D Euler equations with   
finite initial enstrophy.

\section{Two notions of enstrophy defect and Eyink's conjecture}  
  
In this section we will introduce two notions of enstrophy  
defect, one associated with enstrophy dissipation due to viscosity and another associated  
with enstrophy disappearance due to irregular transport. We will also state precisely a version  of Eyink's conjecture
in the setting of full-plane flow.    
  
Let $\omega \in L^{\infty}([0,T); L^{4/3}(\real^2))$ be a weak   
solution of \eqref{e:euler}. Set $j_{\epsilon}(x)=\epsilon^{-2}j(\epsilon^{-1}x)$ to be a   
Friedrichs mollifier and write   
\[  \omega_{\epsilon}=j_{\epsilon}\ast \omega,\]  
\[  u_{\epsilon}=j_{\epsilon}\ast u,\]  
\[  (u \, \omega)_{\epsilon}=j_{\epsilon}\ast (u \, \omega).\]  
  
Then $\omega_{\epsilon}$ solves   
  
\begin{equation} \label{e:mollif1}  
 \begin{aligned}  
  &\partial_{t} \omega_{\epsilon} + \dive \left[ u_{\epsilon}   
   \omega_{\epsilon} + \left((u\omega)_{\epsilon}-u_{\epsilon}  
   \omega_{\epsilon}\right)\right]=0, \\  
  &\omega_{\epsilon}(0) = j_{\epsilon} \ast \omega_0.\\  
 \end{aligned}  
\end{equation}  
  
The associated enstrophy density $\oOmega_{\epsilon}(x,t) = |\omega_{\epsilon}(x,t)|^2 /2 $ satisfies  
  
\begin{equation} \label{e:mollif2}  
  \partial_{t} \oOmega_{\epsilon} + \dive\left[ u_{\epsilon}   
   \oOmega_{\epsilon} + \omega_{\epsilon}\left( (u\omega)_{\epsilon}-u_{\epsilon}  
   \omega_{\epsilon}\right)\right] = -Z_{\epsilon}(\omega),  
\end{equation}  
  
where   
\[  
  Z_{\epsilon}(\omega)= -\nabla\omega_{\epsilon}\cdot  \left(   
  (u\omega)_{\epsilon}-u_{\epsilon} \omega_{\epsilon}\right).  
\]  
The behavior of $Z_{\epsilon}$ as $\epsilon \to 0$ is a description of the space-time distribution of enstrophy   
dissipation of the weak solution $\omega$ due to irregular
transport. We use this notion to define the enstrophy defect.   
  
\begin{definition} \label{ensdefect} The transport enstrophy defect associated to $\omega$ is:  
\[Z^T(\omega) \equiv \lim_{\epsilon \to 0} Z_{\epsilon}(\omega),\]  
whenever the limit exists in the sense of distributions.  
The weak solution $\omega$ is said to be {\em dissipative}  
if  $Z^T(\omega)$ exists and $Z^T(\omega) \geq 0$.  
\end{definition}  
   
Given that the transport enstrophy defect is intended to describe the space-time structure of enstrophy  
dissipation and taking into account that finite-enstrophy weak solutions conserve enstrophy,  one would hope 
that $Z^T(\omega) \equiv 0$ if $\omega_0 \in L^2_c(\real^2)$. Actually, this seems to be a difficult 
problem, to which we will return later on in this work. Recall that, in the 3D case, it is known  
that finite energy solutions may dissipate energy, see \cite{DR,Shni}.   
  
From a physical point of view it is natural to consider weak solutions arising through the vanishing viscosity  
limit. We denote by $\omega_{\nu}$ the solution to the  
two-dimensional Navier-Stokes equations in velocity-vorticity form:  
\begin{subequations} \label{e:visc1}  
 \begin{align}  
  & \partial_{t} \omega_{\nu} + u_{\nu}\cdot \nabla \omega_{\nu} -  
  \nu \Delta \omega_{\nu}=0, \label{e:visc1.a}\\  
  & u_{\nu}=K\ast \omega_{\nu}, \label{e:visc.b}  
 \end{align}  
\end{subequations}  
with initial data $\omega_{0}$. Note that $\dive u_{\nu} = 0$.

The Navier-Stokes evolution naturally dissipates enstrophy, though only through diffusion. 
The  viscous enstrophy density $\oOmega_{\nu}$ satisfies the following parabolic equation: 
\begin{equation} \label{e:visc2}  
  \partial_{t} \oOmega_{\nu} + u_{\nu}\cdot \nabla \oOmega_{\nu} -  
  \nu \Delta \oOmega_{\nu}=-Z^{\nu}(\omega_{\nu}),  
\end{equation}  
where 
\[  
   Z^{\nu}(\omega_{\nu})=\nu |\nabla \omega_{\nu}|^{2}.  
\]  
Note that $Z^{\nu}(\omega_{\nu})\geq 0$ always. We use $Z^{\nu}$ to define a viscous enstrophy 
defect. Let $\omega = \omega(x,t) \in L^{\infty}([0,T);L^p(\real^2))$, $p \geq 4/3$, be a weak solution of  
the 2D Euler equations which was obtained as a vanishing viscosity limit.  
More precisely, we assume that $\omega$ is a limit of a sequence of solutions to the 2D Navier-Stokes equations  
\eqref{e:visc1} with fixed initial data $\omega_0$ and with viscosity $\nu_k \to 0$. In what follows we  
will refer to such a weak solution as a {\em viscosity solution}. Let $\omega_{\nu_k}$ be such an  
approximating sequence of solutions, with   
$\omega_{\nu_k} \rightharpoonup \omega$, weak-$\ast$ in $L^{\infty}([0,T);L^p(\real^2))$.  
Henceforth we will abuse terminology and identify the sequence  
$\{\omega_{\nu_k}\}$ with its weak (inviscid) limit $\omega$.  
 
\begin{definition} \label{vensdefect}  
The viscous enstrophy defect associated to $\omega$ is defined as:  
\[Z^{V}(\omega) \equiv \lim_{\nu_k \to 0} Z^{\nu_k}(\omega_{\nu_k}),\]  
whenever the limit exists in the sense of distributions.  
\end{definition}   
 
Before we formulate Eyink's conjecture we show that, if the initial vorticity has finite enstrophy, 
then the viscous enstrophy defect vanishes identically.

\begin{proposition} \label{viscl2}  
Let $\omega_{0}\in L^{2}_c(\real^{2})$. Let $\omega \in L^{\infty}([0,T);L^2(\real^2))$ be a  
viscosity solution with initial vorticity $\omega_0$. Then $Z^V(\omega)$ exists and it is identically zero.  
\end{proposition}  
  
\begin{proof} Suppose that the viscosity solution $\omega = \omega(x,t)$  
is the limit of the approximating sequence $\omega_{\nu_k}$ of  
solutions to the Navier-Stokes equations. We may assume that $\omega_{\nu_k} \to \omega$  
in $C([0,T), w-\ltwo)$, where $w-\ltwo$ is $L^2$ endowed with the weak topology, see  
\cite{Lio}, Appendix C.  
Multiplying \eqref{e:visc1} by $\omega_{\nu_k}$, integrating by parts, and   
using the divergence-free condition on $u_{\nu_k}$, gives for each fixed  
$\nu_k$ and $t>0$  
\[    \frac{d}{dt} \int_{\real^{2}} \omega_{\nu_k}^{2}(t) \,dx + 2\nu_k  
   \int_{\real^{2}} |\nabla \omega_{\nu_k}|^{2}(t)\,dx = 0. \] 
By integrating in time, we then obtain the same energy estimate as for  
the heat equation, namely  
\begin{equation} \label{e:l2def2}  
  \|\omega_{\nu_k}(t)\|_{\ltwo}^{2} -\|\omega_{0}\|_{\ltwo}^{2} = - 2\nu_k  
  \int_{0}^{t} \int_{\real^{2}}  |\nabla \omega_{\nu_k} |^{2}\,dx\,ds,  
  \quad \forall\, 0<t<T.  
\end{equation}  
 
From Proposition \ref{weakrenorm} it follows that $\omega$ is a renormalized solution  
to \eqref{e:euler.a} and hence $\|\omega(t)\|_{L^2}^2 = \|\omega_0\|_{L^2}^2$. 
Therefore, if $\omega_{\nu_k}(t)$ converges {\em strongly} in  
$\ltwo$ to $\omega(t)$, for each $0<t<T$, then we have that   
\begin{equation} \label{Zkto0}  
   \lim_{\nu_k\to 0}  \int_{0}^{t} \int_{\real^{2}}    
 \nu_k  |\nabla \omega_{\nu_k}|^{2}\,dx\,ds=0.  
\end{equation} 
This means in particular that $\lim_{\nu_k\to 0} Z^{\nu_k}(\omega_{\nu_k})= Z^V(\omega) \equiv 0$ in the  
sense of distributions.  
  
To establish strong convergence of the approximating sequence, we notice that, from \eqref{e:l2def2}, 
$\|\omega_{\nu_k}(t)\|_{\ltwo}\leq \|\omega_{0}\|_{\ltwo}$ for each $t>0$, so that   
\[  
   \limsup_{\nu_k \to 0} \|\omega_{\nu_k}(t)\|_{\ltwo} \leq   
   \|\omega_{0}\|_{\ltwo}=\|\omega(t)\|_{\ltwo}.  
\]  
On the other hand, it follows from the weak lower semicontinuity of the norm that  
\[ \liminf_{\nu_k \to 0} \|\omega_{\nu_k}(t)\|_{\ltwo} \geq   
   \|\omega(t)\|_{\ltwo}, \] 
as $\omega_{\nu_k} \to \omega$ in $C([0,T);w-L^2)$. Thus  
$\|\omega_{\nu_k}(t)\|_{L^2} \to \|\omega(t)\|_{L^2}$ for each $0<t<T$, from which the desired strong  
convergence follows. 
  
\end{proof}  
  
\begin{remark} In the case of periodic flow it is possible to show that, if $\omega$ is a 
dissipative weak solution in $L^{\infty}([0,T);L^2)$, then the transport 
enstrophy defect $Z^T(\omega)$ vanishes identically. The proof is an 
easy adaptation of what was presented above. For the full plane, there are serious technical 
difficulties with controlling the behavior of $Z_{\epsilon}$ near infinity, which are connected 
with understanding the possibility of enstrophy leaving the compact parts of the plane. 
The main concern of the present article is with local enstrophy dissipation so we will avoid this  
issue of escape to infinity. 
\end{remark} 

\vspace{.5cm} 
 
Turbulence theory requires flows that dissipate enstrophy at a rate which does not vanish 
as viscosity goes to zero. A vanishing viscous enstrophy defect excludes precisely such flows. 
From Proposition \ref{viscl2}, we see that in order to model two-dimensional turbulence, one 
should consider flows with infinite enstrophy. Is it possible for flows with infinite enstrophy  
to dissipate enstrophy in a meaningful way? This is the main point in Eyink's work and 
it is precisely what we wish to explore.  
 
A natural choice of space which allows for infinite enstrophy is the $L^2$-based Besov space  
$B^0_{2,\infty}$. The choice of the space $B^{0}_{2,\infty}$ is motivated by the   
Kraichnan-Batchelor theory of two-dimensional turbulence, which  
predicts, in the limit of vanishing viscosity,  
an energy spectrum of the form  
\begin{equation} \label{e:KB}  
  E(\kappa,t)\sim \eta(t)^{2/3}\,\kappa^{-3}.  
\end{equation}  
Above, $\eta(t)$ is the average rate of enstrophy dissipation per unit volume, and   
$E(\kappa,t)$ is the density of the measure $\mu$ given by  
\[  
    \mu(A)=\int_{A} E(\kappa,t)\,d\kappa=\int_{A\times S^{1}} |\Hat{u}(k,t)|^{2}\,dk,  
\]  
with $\kappa=|k|$, for any measurable subset $A$ of the real line.

By Calderon-Zygmund, $\omega \in L^{2}([0,T];B^{0}_{2,\infty})$   
implies that the velocity $u\in$ \newline   
$L^{2}([0,T]; B^{1}_{2,\infty})$ locally, and \cite{Trieb}  
\begin{equation}  
   \|u\|_{B^{1}_{2,\infty}}^{2} \approx \sup_{0<s\leq 1}  
     s^{2} \|\psi(s \cdot) \Hat{u}\|^{2}_{L^{2}},  
\end{equation}  
for $\psi$ a smooth cut-off function supported in the dyadic shell    
$\{k \mid 1/2<|k|<2\}$. By rescaling (here $s=\kappa^{-1}$),   
a finite $B^{1}_{2,\infty}$ norm gives a decay rate like \eqref{e:KB}  
for the energy spectrum.  
  
In this situation, Eyink's conjecture embodies the expectation that the transport enstrophy 
defect accounts for the residual rate of viscous enstrophy dissipation
in the limit of vanishing viscosity.  
One of the main results in the present work is an example showing that this is not necessarily the case.  
  
\begin{conjecture}[Eyink]  
 Let  $\omega$ be a weak solution of the incompressible 2D Euler equations, 
obtained by the vanishing viscosity method, such that $\omega \in L^2((0,T);B^0_{2,\infty}(\real^2))$. 
We assume that there exists $\omega_{\nu_k}$, solutions of the incompressible Navier-Stokes equations  
\eqref{e:visc1}, such that 
\[\omega_{\nu_k} \rightharpoonup \omega \mbox{ in weak-}\ast \;L^2((0,T);B^0_{2,\infty}(\real^2))\]   
Then both limits, $\;\lim_{\nu\to 0^{+}} Z^{\nu}(\omega_{\nu})$ and  
$\lim_{\epsilon\to 0^{+}} Z_{\epsilon}(\omega)$, exist and are equal, so that we may write  
$Z(\omega)=Z^V(\omega)=Z^T(\omega)$ in this case. 
Furthermore, $\omega$ is a dissipative solution. Lastly, there exist one such $\omega$  
with $Z(\omega)>0$.  
\end{conjecture}  

The space $B^0_{2,\infty}$ has the disadvantage of not being 
rearrangement-invariant, which means that it provides no natural estimate for vorticity. 
In addition, $B^0_{2,\infty}$ is not contained in $L^{4/3}$, so that a weak  
solution in this Besov space has to be defined in a different way than what we 
did in Definition \ref{weakvort}, namely using the weak velocity formulation as in \cite{DiPM1}. 

From an analytical standpoint, it is natural to reformulate Eyink's conjecture replacing $B^0_{2,\infty}$
by a rearrangement invariant space containing $L^2$. In that case, the existence of a viscosity weak solution 
follows from appropriate hypotheses on initial data, so that the statement of the
conjecture would become simpler. One straighforward choice is the Marcinkiewicz space 
$L^{2,\infty}$, which is rearrangement invariant. Additionally, vorticities in $L^{2,\infty}$ which are  
supported in sets of finite measure also belong to $L^{4/3}$, so that Definition \ref{weakvort} can be used. 
Although $L^{2,\infty}$ and $B^{0}_{2,\infty}$ are both endpoints of secondary scales of spaces based on $L^2$, 
the precise relation between them has not been clearly stated in the literature. 

The conjecture stated above differs from Eyink's original formulation in that it  
refers to full-plane instead of periodic flow, a distinction which is more technical than substantive. 
One of the purposes of the present article is to produce an example of a weak solution, under the 
constraints of the conjecture, for which both $Z^T$ and $Z^V$ exist, $Z^T \equiv 0$ but $Z^{V}$ does not 
vanish. The example we will present belongs to $L^{2,\infty} \cap B^{0}_{2,\infty}$. 
Before we present the construction of this example, we will examine in more detail the behavior  
of the enstrophy defects in the case of finite enstrophy. This is the subject of the next two sections.

\section{Transport enstrophy defect and local balance of enstrophy} 
 
We have established that, if the initial vorticity has finite enstrophy, then the 
(renormalized) weak solution conserves enstrophy exactly (Proposition \ref{weakrenorm}  
and subsequent observation) and that for viscosity solutions,  the viscous enstrophy defect $Z^V$ vanishes.  
This result implies that, for modeling 2D turbulence, flows with bounded enstrophy are not useful, since they cannot support a cascade. However, independently from its physical relevance, the idea of transport enstrophy defect is very intriguing from the point of view of nonlinear PDE. One of the most interesting problems is whether transport enstrophy dissipation occurs at all, a nontrivial open question .  
In \cite{Ey}, Eyink proved that if the vorticity is in $L^p$, $p>2$ then $Z^T$ exists and vanishes identically. Our main purpose in this section is to examine transport enstrophy dissipation in more detail, looking for criticality in spaces which are logarithmic perturbations of $L^2$.

We begin by considering local balance of enstrophy. One of the ways in
which this balance can be expressed is by showing  
that the enstrophy density $\oOmega$ satisfies the transport equation $\oOmega_t + u \cdot \nabla \oOmega = 0$. 
We first note that, if the initial vorticity $\omega_0$ belongs to $L^2_c$ and if $\omega$ is any weak  
solution with initial vorticity $\omega_0$, then the corresponding enstrophy density $\oOmega = |\omega|^2/2$  
is a {\it renormalized}  
solution of the above transport equation. The proof of this fact follows from the knowledge that  $\omega$ itself is a  
renormalized solution (in this case) and that, if $\beta(s)$ is an admissible renormalization, then so is  
$\beta(s^2)$. This observation is a Lagrangian perspective on local enstrophy balance, but it cannot be immediately translated 
into Eulerian information. We cannot prove that $\oOmega$ is a weak (distributional) solution of the same transport 
equation because of the difficulty in making sense of the term $u \oOmega$ for arbitrary $L^2$ vorticity. 
This difficulty arises since, if the vorticity is in $L^2$, then the associated velocity is only $H^1_{\loc}$ and hence not necessarily bounded.  We will explore this issue further in the following section through examples. 
Our next result is an attempt to  determine the critical space in which viscosity solutions have  
enstrophy densities that solve the transport equation in the sense of distributions. The key idea is to 
identify a critical space where we can make sense of the nonlinear term $u\oOmega$.   
Let us begin by recalling some basic facts regarding  Orlicz and Lorentz spaces. 
 
Let $f \in L^1_c(\real^2)$ and denote by $\lambda_f = \lambda_f(s) \equiv |\{ x \; | \; |f(x)|> s\}|$ its  
distribution function.  
Let $f^{\ast}$ denote the standard nonincreasing rearrangement of $f$,  
see \cite{BS88} for details. We consider the Lorentz spaces $L^{(1,q)}_{\loc}$, based on the maximal function of  
$f^{\ast}$, $f^{\ast\ast}(s) = \frac{1}{s} \int_0^s f^{\ast}(r)dr$, $1\leq q<\infty$: 
\begin{equation} \label{L1q} 
L^{(1,q)}_{\loc} (\real^2) \equiv \{ f \in L^1_c(\real^2) \; | \; \| s \, f^{\ast\ast}(s) \|_{L^q (ds/s)} < \infty \}. 
\end{equation}  
 
There are two ways of defining Lorentz spaces, one based on $f^{\ast\ast}$ and the other based on $f^{\ast}$.  
The two definitions are equivalent if $p>1$, but they lead to two slightly different spaces if $p=1$, which 
are usually denoted $L^{1,q}$ and $L^{(1,q)}$. The spaces $L^{(1,q)}_{\loc}$ play a distinguished role  
in the study of incompressible 2D  Euler: if $ 1 \leq q < 2$  they can be  
compactly imbedded in $H^{-1}_{\loc}$. If $q = 2$ then the imbedding is merely continuous, see \cite{LLT}. In fact,  
it was observed by P.-L. Lions in \cite{Lio} that $L^{(1,2)}_{\loc}(\real^2)$ is the {\em largest} rearrangement  
invariant Banach space that can be continuously imbedded in $H^{-1}_{\loc}(\real^2)$.  
 
Let $1\leq p < \infty$ and $a \in \real$. Define $A_{p,a} = A_{p,a}(s) \equiv [s \log^{a} (2+s)]^p$, for $s>0$. Then this  
is a $\Delta$-regular $N$-function (see \cite{Adams} for the basic definitions). In particular $A_{p,a}$ is  
nondecreasing and convex. The associated Orlicz space is the Zygmund space $L^p(\log L)^a$ defined by: 
\begin{equation} \label{LplogLa} 
L^p(\log L)^a (\real^2) \equiv \left\{f \in L^1_{\loc} \, \left| \, \int_{\real^2} A_{p,a}(|f(x)|) \, dx \right. <  
\infty.\right\} 
\end{equation} 
 
The Orlicz spaces are Banach spaces when equipped with the Luxemburg norm: 
\begin{equation} \label{luxnorm} 
\|f\|_{p,a} = \inf \left\{ k>0 \, \left| \, \int A_{p,a}\left( \frac{|f(x)|}{k}\right.\right) \, dx \leq 1\right\}.  
\end{equation} 
If $f$ does not vanish identically then the infimum is attained.  
 
If $p = 1$, these spaces are well-known logarithmic refinements of $L^1$ commonly denoted by $L(\log L)^a$;  
for arbitrary $p$ they are logarithmic refinements of $L^p$. It was observed in  
\cite{LLT} that $L(\log L)^{1/q} \subset L^{(1,q)} \subset L(\log L)^a$ for any $a < 1/q \leq 1$.  
The relevant case at present is $q=2$.  
 
We begin with a technical lemma. 
 
\begin{lemma} \label{product} 
Let $\alpha$ and $\beta$ be functions in $L^2(\log L)^{1/4}(\real^2)$. Then the product $\alpha \beta$ belongs to  
$L (\log L)^{1/2}$ and 
\[ \|\alpha \beta \|_{1,1/2} \leq 4 \left(\max \{\|\alpha\|_{2,1/4};\|\beta\|_{2,1/4}\}\right)^2.\] 
\end{lemma} 
 
\begin{proof} 
We may assume without loss of generality that neither $\alpha$ nor $\beta$ vanish identically, otherwise the result is immediate. Thus the infimum in the Luxemburg norm \eqref{luxnorm} is attained for both $\alpha$ and $\beta$, i.e.,  
 
\[\int A_{2,1/4}\left( \frac{|\alpha(x)|}{\|\alpha\|_{2,1/4}}\right) \, dx = 1  \mbox{ and }  
\int A_{2,1/4}\left( \frac{|\beta(x)|}{\|\beta\|_{2,1/4}}\right) \, dx = 1.\] 
It is an easy exercise to show that $A_{2,1/4}(2s) \geq 4 A_{2,1/4}(s)$, for any $s>0$. Thus it follows that  
 
\[\int A_{2,1/4}\left( \frac{|\alpha(x)|}{2\|\alpha\|_{2,1/4}}\right) \, dx \leq \frac{1}{4}  \mbox{ and }  
\int A_{2,1/4}\left( \frac{|\beta(x)|}{2\|\beta\|_{2,1/4}}\right) \, dx \leq \frac{1}{4}.\] 
 
Let $k = \max\{\|\alpha\|_{2,1/4};\|\beta\|_{2,1/4}\}$. Then: 
\[\int A_{1,1/2}\left(\frac{\alpha(x)\beta(x)}{4k^2} \right)\, dx = \] 
\[=\int_{\alpha \geq \beta} A_{1,1/2}\left(\frac{\alpha(x)\beta(x)}{4k^2} \right)\, dx +  
\int_{\beta > \alpha}A_{1,1/2}\left(\frac{\alpha(x)\beta(x)}{4k^2} \right)\, dx\] 
\[\leq \int \frac{|\alpha|^2}{(2k)^2}\log^{1/2}\left( 2 + \frac{|\alpha|^2}{(2k)^2}  \right) \, dx +  
\int \frac{|\beta|^2}{(2k)^2}\log^{1/2}\left( 2 + \frac{|\beta|^2}{(2k)^2}  \right) \, dx \] 
\[\leq \sqrt{2} \int A_{2,1/4}\left( \frac{\alpha}{2k}\right) \, dx +  
\sqrt{2} \int A_{2,1/4}\left( \frac{\beta}{2k}\right) \, dx \] 
\[ \leq \sqrt{2}(1/4 + 1/4) < 1,\] 
where the last estimate holds in view of the fact that $A_{2,1/4}$ is nondecreasing. 
 
It follows that  
\[\|\alpha\beta\|_{1,1/2} \leq 4k^2,\] 
as we wished. 
 
\end{proof} 
 
We are now ready to prove that the enstrophy density is a weak solution of the appropriate transport equation, if the vorticity is an inviscid limit and belongs to $L^2(\log L)^{1/4}$. 
 
\begin{theorem} \label{ensbalance} 
Let $\omega_0 \in (L^2(\log L)^{1/4})_c(\real^2)$. Consider a viscosity solution  
$\omega \in L^{\infty}([0,T);L^2(\log L)^{1/4}(\real^2))$ with initial data $\omega_0$.  
Then the following equation holds in the sense of distributions: 
\begin{equation} \label{e:ensbalance}
 \partial_t (|\omega|^2) + \dive ( u |\omega|^2) = 0,
\end{equation}
where $u = K \ast \omega$. 
\end{theorem} 
  
\begin{proof} 
Let $\omega_{\nu_k}$, $\nu_k \to 0$, be a sequence of solutions of the 2D Navier-Stokes equations \eqref{e:visc1}, with  
initial vorticity $\omega_0$, such that $\omega_{\nu_k} \rightharpoonup \omega$ weak-$\ast$ in  
$L^{\infty}([0,T);L^2(\real^2))$. The existence of such a sequence is guaranteed by the fact that $\omega$ is 
a viscosity solution with initial vorticity $\omega_0 \in L^2(\log L)^{1/4} \subset L^2$.  
 
We will begin by showing an {\it a priori} bound, uniform in viscosity, in the space  
$L^{\infty}([0,T);L^2(\log L)^{1/4}(\real^2))$ for $\omega_{\nu_k}$. 
To this end we multiply \eqref{e:visc1.a} by  
\[\frac{1}{m}A_{2,1/4}^{\prime}\left(\frac{\omega_{\nu_k}}{m}\right),\] 
for arbitrary $m>0$. Here, $A_{2,1/4}^{\prime}$ is the derivative of $A_{2,1/4}$ with respect to its argument.
Then $\frac{1}{m}\omega_{\nu_k}$ satisfies the following equation: 
\begin{equation} \label{a214eq} 
\begin{aligned} 
\partial_t\left(A_{2,1/4}\left(\frac{\omega_{\nu_k}}{m}\right)\right) & + u_{\nu_k} \cdot \nabla  
A_{2,1/4}\left(\frac{\omega_{\nu_k}}{m}\right) = \\ 
& \nu_k \Delta A_{2,1/4}\left(\frac{\omega_{\nu_k}}{m}\right) -  
\frac{\nu_k}{m^2} A_{2,1/4}^{\prime\prime}\left(\frac{\omega_{\nu_k}}{m}\right)|\nabla \omega_{\nu_k}|^2. 
\end{aligned} 
\end{equation} 
We integrate \eqref{a214eq} in all of $\real^2$, use the divergence-free condition on velocity and the convexity of 
$A_{2,1/4}$ to conclude that, for any $m>0$, 
\[\frac{d}{dt}\int_{\real^2}A_{2,1/4}\left(\frac{\omega_{\nu_k}(x,t)}{m}\right)\,dx \leq 0.\] 
Thus, since the norm in $L^2(\log L)^{1/4}$ is the Luxemburg norm \eqref{luxnorm}, it follows that  
\begin{equation} \label{ap23}\|\omega_{\nu_k}(\cdot,t)\|_{2,1/4} \leq \|\omega_0\|_{2,1/4}, 
\end{equation} 
for any $0\leq t<T$.  
 
We have obtained that $\omega_{\nu_k}$ is  bounded in $L^{\infty}([0,T);L^2(\log L)^{1/4}(\real^2))$ and,  
as this is a Banach space, we may assume, passing to a subsequence if necessary, that $\omega_{\nu_k} \rightharpoonup  
\omega$ weak-$\ast$ in this space as $\nu_k \to 0$.  
 
Next recall that $\oOmega_{\nu_k} = |\omega_{\nu_k}|^2/2$ satisfies the viscous enstrophy balance equation \eqref{e:visc2}. Therefore, 
for any test function $\varphi \in \mathcal{D}((0,T)\times\real^2)$ we have: 
 
\begin{equation} \label{vwkform}  
\begin{aligned} 
\int_0^T \int_{\real^2} &\varphi_t \oOmega_{\nu_k} \,dxdt +  
\int_0^T \int_{\real^2} \nabla \varphi \cdot u_{\nu_k} \oOmega_{\nu_k}  \, dxdt= \\ 
\\ 
&=\int_0^T \int_{\real^2} \nu_k \Delta \varphi \, \oOmega_{\nu_k} \,dxdt - 
\int_0^T \int_{\real^2} \varphi Z^{\nu_k}(\omega_{\nu_k}) \, dxdt.\\ 
\end{aligned} 
\end{equation} 
 
We need to pass to the limit $\nu_k \to 0$ in each of the terms above. First recall,  
from the proof of Proposition \ref{viscl2}, that $\oOmega_{\nu_k}(\cdot,t) \to \oOmega(\cdot,t)$ strongly in  
$L^1(\real^2)$ for each $0<t<T$. Indeed, we used this fact to show that $Z^{\nu_k}(\omega_{\nu_k}) \to 0$ 
in $L^1([0,T)\times\real^2)$, see \eqref{Zkto0}. Furthermore, as $\int\oOmega_{\nu_k}(\cdot,t)\,dx \leq \int\oOmega_0\,dx$ it follows, by the Dominated Convergence Theorem,  
that $\oOmega_{\nu_k} \to \oOmega$ strongly in $L^1([0,T)\times\real^2)$. Therefore, the 
first term in \eqref{vwkform} converges to 
\[\int_0^T \int_{\real^2} \varphi_t \oOmega \, dxdt, \] 
and the third term converges to zero due to the vanishing factor $\nu_k$. The fourth term in \eqref{vwkform} converges to zero, as was shown in \eqref{Zkto0} in the proof of Proposition \ref{viscl2}. It remains to determine the limit behavior of the nonlinear term. 
 
We start with the observation that $\omega_0 \in (L^2(\log L)^{1/4})_c \subset L^1_c$. Using the maximum principle it is easy to show that the $L^1$-norm of the solution $\omega_{\nu_k}$ decreases in time: 
\begin{equation} \label{L1estimate} 
 \|\omega_{\nu_k}(\cdot, t)\|_{L^1(\real^2)} \leq \|\omega_0\|_{L^1(\real^2)}. 
\end{equation} 
Thus, as the Biot-Savart kernel $K$ is locally integrable and bounded near infinity, the convolution  
$K\ast \omega_{\nu_k}$ is well-defined. 
We may therefore use the Biot-Savart law $u_{\nu_k} = K \ast \omega_{\nu_k}$ to find: 
\begin{equation} \label{nlterm} 
\int_0^T \int_{\real^2}    \nabla \varphi \cdot u_{\nu_k} \oOmega_{\nu_k}  \, dxdt =  
\end{equation} 
\[ 
\begin{aligned} 
& =\int_0^T \int_{\real^2}   \int_{\real^2}  \nabla \varphi(x,t) \cdot K(x-y) \omega_{\nu_k}(y,t) \oOmega_{\nu_k}(x,t)  \, dydxdt \\ 
& =-\int_0^T  \int_{\real^2} \omega_{\nu_k}(y,t)  \int_{\real^2} K(y-x)\cdot\nabla\varphi(x,t)\oOmega_{\nu_k}(x,t)\,dx  \,  
dydt,  
\end{aligned} 
\]  
as $K$ is antisymmetric. Thus we may write  
\[ \eqref{nlterm} = - \int_0^T \int_{\real^2} \omega_{\nu_k}(y,t) \,{\mathcal I}_k(y,t)\,dydt, \] 
with  
\[ {\mathcal I}_k \equiv  
\int_{\real^2}  
K(y-x)\cdot\nabla\varphi(x,t)\oOmega_{\nu_k}(x,t)\,dx. \] 
 
Let $\oOmega = |\omega|^2/2$. Denote by ${\mathcal I}$ the function 
\[ {\mathcal I} \equiv  \int_{\real^2}  
K(y-x)\cdot\nabla\varphi(x,t)\oOmega(x,t)\,dx, \] 
which is well defined, as we will see later. 
 
We deduce, from the {\it a priori} estimate \eqref{ap23} in $L^2(\log L)^{1/4}$, from Lemma \ref{product}, and from the fact that 
each component of $\nabla \varphi$ is a smooth test function, that $\{\nabla \phi \oOmega_{\nu_k}\}$ is bounded in  
$L^{\infty}((0,T);L(\log L)^{1/2}(\real^2))$ 
and, therefore, in $L^{\infty}((0,T);L^{(1,2)}_{\loc}(\real^2))$ (see 
\cite{LLT}). As already observed above, $L^{(1,2)}_{\loc}$ can be continuously imbedded in $H^{-1}_{\loc}$, so that 
\begin{equation} \label{IkL2loc} 
\{\mathcal{I}_k \} \mbox{ is bounded in } L^{\infty}((0,T);L^2_{\loc}(\real^2)). 
\end{equation}  
 
Thus it follows that, passing to a subsequence if necessary, ${\mathcal I}_k$ converges weak-$\ast$ in 
$L^{\infty}((0,T);L^2_{\loc}(\real^2))$ to a weak limit. We will show that this weak limit is  
${\mathcal I}$. We know that $\nabla \varphi \oOmega_{\nu_k} \to \nabla \varphi \oOmega$ strongly in  
$L^1((0,T)\times \real^2)$, because $\omega_{\nu_k} \to \omega$ strongly in $L^2((0,T)\times\real^2)$.  
Let $\eta \in {\cal D}((0,T)\times\real^2)$. Then we may write: 
\[\langle {\mathcal I}_k , \eta \rangle = -\langle \nabla \varphi \oOmega_{\nu_k}, K \ast \eta \rangle, \] 
using the antisymmetry of $K$. Since $K\ast \eta \in L^{\infty}((0,T)\times\real^2)$ we therefore obtain that  
\[ \langle {\mathcal I}_k , \eta \rangle \to -\langle \nabla \varphi \oOmega , K \ast \eta \rangle =  
\langle {\mathcal I} , \eta \rangle.\] 
We have shown that ${\mathcal I}_k \to {\mathcal I}$ in the sense of distributions, so that, by uniqueness of limits,   
the weak limit of ${\mathcal I}_k$ is necessarily equal to ${\mathcal
  I}$. Hence, the whole sequence ${\mathcal I}_k$ converges weakly to ${\mathcal I}$, without the need to pass to a subsequence. In particular, we have established  
that the integral in the definition of ${\mathcal I}$ is well defined. 
 
The next step is to deal with the behavior of ${\mathcal I}_k$ at infinity.  
Note that each component of $\nabla \varphi \oOmega_{\nu_k}$ is  
compactly supported, uniformly in $t$ and $\nu_k$, in a ball, say,  
$B(0;R)$. As the viscous enstrophy decreases in time, we find that   
\[\|\nabla\varphi\oOmega_{\nu_k}\|_{L^{\infty}((0,T);L^1(\real^2))} \leq  
C(\varphi)  \,\int_{\real^2} \oOmega_0 \, dx \equiv C(\varphi)\,\Omega_0.\] 
From this observation and the explicit expression for the kernel $K$, a direct estimate yields that  
\[ |{\mathcal I}_k(y,t)| \leq \frac{\widetilde{C}(\varphi)\,\Omega_0}{|y|}\] 
for $|y|\geq 2R$. Hence  
\begin{equation} \label{IkLinfty} 
\{{\mathcal I}_k \} \mbox{ is bounded in }  
L^{\infty}((0,T)\times (\real^2 \setminus B(0;2R))). 
\end{equation} 
 
Using the same argument as was used above to establish that ${\mathcal I}_k \rightharpoonup  {\mathcal I}$ weak-$\ast$ 
$L^{\infty}((0,T);L^2_{\loc}(\real^2))$, we may conclude, from estimate \eqref{IkLinfty}, that  
${\mathcal I}_k \rightharpoonup {\mathcal I}$ weak-$\ast$ in $L^{\infty}((0,T)\times (\real^2 \setminus B(0;2R)))$  
as well, without the need to pass to a subsequence.  
 
We claim that $\omega_{\nu_k} \to \omega $ strongly in $L^1((0,T)\times\real^2)$, 
as $\nu_k \to 0$, as well as in $L^2((0,T)\times\real^2)$.  
Assuming the claim, we can pass to the limit in the nonlinear term. Indeed, we write 
\[\eqref{nlterm} = -\left( \int_0^T \int_{B(0;2R)} \omega_{\nu_k}(y,t) \,{\mathcal I}_k(y,t)\,dydt + \right.\]
\[ \left. + \int_0^T \int_{\real^2 \setminus B(0;2R)}\omega_{\nu_k}(y,t) \,{\mathcal I}_k(y,t)\,dydt \right),\] 
which converges to 
\[-\int_0^T \int_{B(0;2R)} \omega(y,t) \,{\mathcal I}(y,t)\,dydt 
- \int_0^T \int_{\real^2 \setminus B(0;2R)}\omega(y,t) \,{\mathcal I}(y,t)\,dydt,\] 
as each integral forms a ``weak-strong pair'', by virtue of the convergence ${\mathcal I}_k \to {\mathcal I}$  
established above, and noting that $L^{\infty}((0,T);L^2_{\loc}) \subset L^2((0,T);L^2_{\loc})$.

All that remains is to prove the claim. We begin by noting that the
strong convergence in $L^2((0,T)\times\real^2)$ was observed in the
proof of Proposition \ref{viscl2}: it follows from the convergence of
the norms together with weak convergence. To address strong
convergence in $L^1$ we make use of the following fact (for $p=1$),
due to H. Br\'{e}zis and E. Lieb, (see Theorem 8 of \cite{EvansCBMS}
for a proof): a sequence that converges weakly and almost everywhere
and such that the $L^p$-norms also converge will converge strongly in
$L^p$. We obtain weak convergence in $L^1((0,T)\times \real^2)$,
passing to a subsequence if necessary, directly from the {\it a
  priori} estimate \eqref{L1estimate} on the $L^1$-norm of $\omega_{\nu_k}$  together with strong convergence in $L^2$. We also have almost everywhere convergence passing to a further subsequence if needed. Finally, we can establish strong convergence of the $L^1$-norm by repeating the argument used in the proof of Proposition \ref{viscl2} to show that the $L^2$-norms converge. Consequently, strong convergence in $L^1$ holds for this particular subsequence. However, since we have identified the limit, we find that the whole sequence $\omega_{\nu_k}$ converges to $\omega$ strongly in $L^1((0,T) \times \real^2)$ as  $\nu_k \to 0$, as we wished. 
 
\end{proof} 
 
\begin{remark} The natural condition under which the argument above remains valid is $|\omega_0|^2 \in L^{(1,2)}$. 
We chose to present the result under the slightly stronger assumption $\omega_0 \in L^2(\log L)^{1/4}$ because 
it does not seem immediate to provide an {\it a priori} estimate on the square of $\omega$ in $L^{(1,2)}$ that is  
uniform in viscosity.  
\end{remark} 

\begin{remark} If we do not assume that the weak solution is a viscosity solution then the best result available 
on enstrophy density satisfying the transport equation
\eqref{e:ensbalance} in the sense of distributions requires $\omega_0 \in L^p_c$, 
$p>2$, see \cite{Ey}. 
\end{remark}

\vspace{0.5cm} 
 
One is naturally led to ask what knowledge has been gained with Theorem \ref{ensbalance}. Where we previously 
knew that the enstrophy density satisfied the transport equation in the renormalized sense, we now know  
that the equation is satisfied in the sense of distributions. We apply this additional information in the 
proof of our next result.  
 
In the remainder of this section, we are concerned with the conditions under which
$Z^T$ exists and vanishes for finite enstrophy flows. The key point in the proof of Theorem \ref{ensbalance} is that we provided meaning to the term $u \oOmega$ for $\omega \in L^2(\log L)^{1/4}$, through the computation of   
\eqref{nlterm}.  Assigning meaning to the nonlinearity $u \oOmega$ will also play a central role in the proof of the next result. We formalize the meaning we wish to adopt in a definition. 
 
\begin{definition} \label{uoomega} 
Let $\omega \in L^2(\log L)^{1/4}(\real^2)\cap L^1(\real^2)$. Let $u = K \ast \omega$. Then 
we define $u \oOmega \in {\mathcal D}^{\prime}(\real^2)$ by:  
\[\langle u \oOmega, \Phi \rangle = - \int_{\real^2} \omega(y) \int_{\real^2} K(y-x) \cdot \Phi(x)\oOmega(x) \,dx \,dy \] 
\[\equiv - \int_{\real^2} \omega(y) [K \cdot \ast (\Phi \oOmega)](y)\, dy, \] 
for any test vector field $\Phi \in {\mathcal D}(\real^2)$. 
\end{definition} 
The integral above is well-defined as $\Phi \oOmega$ is a compactly supported function in  
$L(\log L)^{1/2}_{\loc} \hookrightarrow L^{(1,2)}_{\loc}$ and $\omega
\in L^2 \cap L^1$, see the proof of Theorem
\ref{ensbalance}. Moreover, it is easy to 
establish that $\Phi \mapsto \langle u \oOmega, \Phi \rangle$ is a continuous linear functional over $\mathcal{D}$.
 
We are now ready to state and prove our final result in this section. 
 
\begin{theorem} \label{znl} 
Let $\omega \in L^{\infty}([0,T);L^2(\log L)^{1/4}(\real^2)\cap L^1(\real^2))$ be a weak solution of the incompressible 
2D Euler equations. Then the transport enstrophy defect $Z^T(\omega)$ exists (as a distribution).  
If $\omega$ is a viscosity solution with initial vorticity $\omega_0 \in (L^2(\log L)^{1/4})_c (\real^2)$ then 
 $Z^T(\omega) \equiv 0$.  
\end{theorem} 
 
\begin{proof} 
Let $j_{\epsilon}$ be a radially symmetric, compactly supported Friedrichs mollifier.  
Recall the notation $\omega_{\epsilon}$, $u_{\epsilon}$ and $(u\omega)_{\epsilon}$  
introduced in the beginning of Section 3.  
 
Let $\varphi \in {\cal D}((0,T)\times\real^2)$. We multiply the equation \eqref{e:mollif2} for $\oOmega_{\epsilon}= 
|\omega_{\epsilon}|^2/2$ by $\varphi$ and integrate over $(0,T)\times\real^2$ to find: 
\begin{equation} \label{wkformollens}  
\begin{aligned} 
\int_0^T & \int_{\real^2}\varphi_t \oOmega_{\epsilon} \,dxdt  + \int_0^T\int_{\real^2}\nabla\varphi \cdot u_{\epsilon}   
   \oOmega_{\epsilon} \, dxdt +\\ 
&+ \int_0^T\int_{\real^2} \nabla\varphi \,\omega_{\epsilon}\cdot \left( (u\omega)_{\epsilon}-u_{\epsilon} \omega_{\epsilon}\right)\, dxdt = \int_0^T\int_{\real^2}\varphi Z_{\epsilon}(\omega)\,dxdt.  
\end{aligned} 
\end{equation}  
We wish to pass to the limit $\epsilon \to 0$.  Let us begin by examining the first two terms above. 
 
The integrand in the first term is $\varphi_t |\omega_{\epsilon}|^2/2$, which converges to $\varphi_t|\omega|^2/2$, as $\epsilon \to 0$, strongly in $L^1((0,T)\times\real^2)$. Indeed, by standard properties of mollifiers, $\omega_{\epsilon}(\cdot,t) \to \omega(\cdot,t)$ strongly in $L^2(\real^2)$ for each $0<t<T$, and also  
\[\|\omega_{\epsilon}(\cdot,t)\|_{L^2}\leq\|\omega(\cdot,t)\|_{L^2} \equiv \|\omega_0\|_{L^2}.\]  
Hence we may obtain the desired conclusion using the Dominated Convergence Theorem.  
 
Next we consider the second term. We note that mollification is continuous in $\Delta$-regular Orlicz spaces (see Theorem 8.20 in \cite{Adams}) so that $\omega_{\epsilon} \in L^{\infty}([0,T);
L^2(\log L)^{1/4}(\real^2)\cap L^1(\real^2))$. As convolutions are associative, we have that  
$u_{\epsilon}=K\ast\omega_{\epsilon}$. We are thus in position to write the second term in \eqref{wkformollens} using  
Definition \ref{uoomega}: 
\begin{equation}\label{bluuug} 
\int_0^T\int_{\real^2}\nabla \varphi \cdot u_{\epsilon}  \oOmega_{\epsilon} \, dxdt =  
- \int_0^T \int_{\real^2} \omega_{\epsilon}(y,t) [K \cdot \ast (\nabla \varphi \oOmega_{\epsilon})](y,t)\, dy dt. 
\end{equation} 
 
It follows from Lemma \ref{product} that the family $\nabla \varphi \oOmega_{\epsilon}$ is uniformly bounded in  
$L^{\infty}((0,T);L(\log L)^{1/2}(\real^2))$. Hence we find, as in  \eqref{IkL2loc}, that  
\[\left\{ K \cdot \ast (\nabla \varphi \oOmega_{\epsilon}) \right\} 
 \mbox{ is bounded in } L^{\infty}((0,T);L^2_{\loc}(\real^2)).\] 
Furthermore, $\|\nabla \varphi \oOmega_{\epsilon}(\cdot,t)\|_{L^1} \leq \|\nabla \varphi \oOmega_0\|_{L^1}$  and 
$\nabla \varphi \oOmega_{\epsilon}$ has compact support {\em uniformly} in $t$ and $\epsilon$, so that, as in  
\eqref{IkLinfty}, 
\[\left\{ K \cdot \ast (\nabla \varphi \oOmega_{\epsilon}) \right\} 
\mbox{ is bounded in } L^{\infty}((0,T)\times(\real^2\setminus B(0;2R))),\] 
for $R$ sufficiently large.
Standard properties of mollifiers yield that $\omega_{\epsilon} \to \omega$ strongly in both   
$L^2((0,T)\times\real^2)$ and $L^1((0,T)\times\real^2)$. Thus we may conclude, as in the proof of  
Theorem \ref{ensbalance}, that the left hand side of \eqref{bluuug} converges to  
\[- \int_0^T \int_{\real^2} \omega(y,t) [K \cdot \ast (\nabla \varphi \oOmega)](y,t)\, dy dt, \] 
when $\epsilon \to 0$.  
 
Finally, let us examine the third term. The key point in this proof is to show that it vanishes as $\epsilon \to 0$.  
We use  
the radial symmetry of the mollifier $j_{\epsilon}$ to obtain: 
\[\int_0^T\int_{\real^2} \nabla\varphi \, \omega_{\epsilon} \cdot  
\left( (u\omega)_{\epsilon}-u_{\epsilon} \omega_{\epsilon}\right)\, dxdt = \] 
\[\int_0^T\int_{\real^2} [(\nabla\varphi\,\omega_{\epsilon})\ast j_{\epsilon} ]\cdot u \omega \, dxdt -  
\int_0^T\int_{\real^2}  \nabla\varphi \cdot u_{\epsilon} |\omega_{\epsilon}|^2\ \, dxdt \equiv {\mathcal I}_{\epsilon}  
- {\mathcal J}_{\epsilon}.\] 
We have already analyzed ${\mathcal J}_{\epsilon}$ in \eqref{bluuug}. We know that 
\[{\mathcal J}_{\epsilon} \equiv  
-\int_0^T\int_{\real^2} \omega_{\epsilon} (y,t)[K\cdot\ast(\nabla\varphi|\omega_{\epsilon}|^2)](y,t)\,dydt \] 
\[\to  
-\int_0^T \int_{\real^2} \omega(y,t) [K \cdot \ast (\nabla \varphi |\omega|^2)](y,t)\, dy dt, \] 
as $\epsilon \to 0$. We will now analyze ${\mathcal I}_{\epsilon}$.  
We start by observing that, using the antisymmetry of $K$, we can write: 
\[{\mathcal I}_{\epsilon} = -\int_0^T\int_{\real^2} \omega(y,t)  
[ \,K\cdot\ast ( \,\omega ( \,( \nabla\varphi\,\omega_{\epsilon}) \ast j_{\epsilon}\, )\, )\,] 
(y,t)\,dydt.\] 
Next we note that, by standard properties of mollification, $(\nabla\varphi\,\omega_{\epsilon})\ast j_{\epsilon} \to  
\nabla \varphi \,\omega$ strongly in $L^2((0,T)\times\real^2)$ as $\epsilon \to 0$. In addition,  
$(\nabla\varphi\,\omega_{\epsilon})\ast j_{\epsilon}$ is compactly supported, uniformly in $t$ and $\epsilon$, and it is  
uniformly bounded in $L^{\infty}((0,T);L^2(\log L)^{1/4}(\real^2))$.  
Using  Lemma \ref{product} we deduce that 
\[\{\omega(\,(\nabla\varphi\,\omega_{\epsilon})\ast j_{\epsilon}\,) \} \mbox{ is bounded in } 
L^{\infty}((0,T);L(\log L)^{1/2}(\real^2)).\] 
Therefore $\{K\cdot\ast ( \,\omega ( \,( \nabla\varphi\,\omega_{\epsilon}) \ast j_{\epsilon}\, )\, ) \}  
\mbox{ is bounded in } 
L^{\infty}((0,T);L^2_{\loc}(\real^2))$ and in  
$L^{\infty}((0,T)\times(\real^2\setminus B(0;2R)))$, for $R$ sufficiently large. 
From this observation we may conclude, as we have before, that  
\[ {\mathcal I}_{\epsilon} \to  
-\int_0^T\int_{\real^2} \omega(y,t) [K \cdot \ast (\nabla\varphi|\omega|^2)](y,t)\,dy dt, \] 
as $\epsilon \to 0$. Therefore the third term vanishes in the limit $\epsilon \to 0$. The proof that  
the transport enstrophy defect exists as a distribution is complete. In fact, we have established that  
\begin{equation} 
\begin{aligned}\label{ztlimit} 
\langle Z^T(\omega), \varphi \rangle & =  
\int_0^T\int_{\real^2}\varphi_t \oOmega \, dxdt  
- \int_0^T \int_{\real^2} \omega(y,t) [K \cdot \ast (\nabla \varphi \oOmega)](y,t)\, dy dt \\  
& = \langle \oOmega , \varphi_t \rangle + \langle u\oOmega , \nabla \varphi \rangle  
\equiv - \langle \oOmega_t + \dive (u\oOmega) , \varphi \rangle,  
\end{aligned} 
\end{equation} 
where the former identity follows from Definition \ref{uoomega}. 
 
Finally, in view of Theorem \ref{ensbalance} we have that, if $\omega$ is a viscosity solution 
with initial vorticity $\omega_0 \in (L^2(\log L)^{1/4})_c(\real^2)$, then the enstrophy density 
balance equation holds in the sense of distributions, so \eqref{ztlimit} above implies $Z^T(\omega) \equiv 0$  
in this case. 
 
\end{proof} 
         
\begin{remark} This result raises a few interesting questions. First, if one could find an example of a weak solution  
with initial vorticity in $L^2(\log L)^{1/4}$ and such that $Z^T$ does not vanish, one would have established  
nonuniqueness of weak solutions. In fact, any example where $Z^T$ exists and does not vanish would be quite interesting.  
Second, one naturally wonders how sharp is the regularity condition $L^2 (\log L)^{1/4}$ on vorticity. This is the  
subject of the next section.  
\end{remark}

\section{The Biot-Savart law in $L^2$-based Zygmund spaces}  
  
The purpose of this section is to illustrate the behavior of the term 
$u|\omega|^2$ through examples. We will be considering pairs $(u,\,\omega)$ 
related by the Biot-Savart law, but not necessarily solutions of the 2D Euler equations. 
We will not establish that the condition $\omega \in L^2(\log L)^{1/4}$ 
(or $|\omega|^2 \in L^{(1,2)}$) is necessary 
for making sense of $u|\omega|^2$, but we will exhibit an example showing
that it is not possible to define $u|\omega|^2$ as a distribution for an arbitrary
vorticity in $L^2$. Furthermore, the family of examples we will present also 
proves that the velocities associated to vorticities in $L^2(\log L)^{1/4}$ are
not necessarily bounded, something which would trivialize the proofs in the 
previous section.

It would be natural to look for such examples in the class of radially 
symmetric vorticities, but we will see in our first Lemma that this approach is
not useful.

\begin{lemma} \label{symbd}
Let $\omega \in L^2_c(\R^2)$, such that $\omega (x) = \phi(|x|)$. 
Let $u \equiv K \ast \omega$. Then $u$ is bounded and 
$\|u\|_{L^{\infty}}\leq C\|\omega\|_{L^2} $.
\end{lemma}  

\begin{proof}
The reader may easily check that if the vorticity is radially symmetric, 
then the Biot-Savart law becomes: 
\[ u(x) = \frac{x^{\perp}}{|x|^2} \int_0^{|x|} s \phi(s) ds. \] 
As $\omega \in L^2$, it follows that $\phi \in L^2(sds)$. We use the 
Cauchy-Schwarz inequality with respect to $sds$ to obtain: 
\[\left|\int_0^{|x|} s \phi(s) ds \right|\leq \left(\int_0^{|x|} s ds \right)^{1/2}  
\left( \int_0^{|x|} s \phi^2(s) ds \right)^{1/2} \leq C \|\omega\|_{L^2} |x|.\] 
This concludes the proof. 
\end{proof}

Recall that the velocity associated to an $L^p$ vorticity is bounded if $p>2$,
but logarithmic singularities may occur when $p=2$. The symmetry in a radial 
vorticity configuration implies a certain cancellation in the Biot-Savart law,
and it is this cancellation  which is responsible for the additional 
regularity observed in the lemma above. We will consider a family of examples
given by breaking the symmetry in the simplest way possible. 

Let $1/2< \alpha <1$. We will denote by $\omega^{\alpha}_+$ the function
\[\omega^{\alpha}_+(x) \equiv \frac{1}{|x| |\log |x||^{\alpha}} 
\chi_{B^+(0;1/3)}(x), \]
where $B^+(0;1/3) = B(0;1/3) \cap \{x_2 > 0\}$.

Note first that $\omega^{\alpha}_+ \in L^2_c$. Indeed,
\[\| \omega^{\alpha}_+ \|_{L^2}^2 = \pi \int_0^{1/3} 
\frac{ds}{s|\log s|^{2\alpha}} = \frac{\pi}{2\alpha -1} (\log 3)^{1-2\alpha},\]
as long as $\alpha > 1/2$. We can make a more precise characterization of
the regularity of $\omega^{\alpha}_+$ using the Zygmund class hierarchy.
 
We denote the radially symmetric extension of $\omega^{\alpha}_+$ as 
\begin{equation} \label{cizania} 
\omega^{\alpha}(x) \equiv (|x||\log |x||^{\alpha})^{-1}\chi_{B(0;1/3)}. 
\end{equation}

\begin{lemma} \label{precreg}
If $1/2 < \alpha < 1$ then $\omega^{\alpha}_+ \in L^2(\log L)^{\kappa}$, for all
$0 \leq \kappa < \alpha - 1/2$. 
\end{lemma}

\begin{proof}
We observe that $(|x||\log |x||^{\alpha})^{-1}$ is a decreasing function of $|x|$ if $|x|\leq e^{-\alpha}$.  
In particular, as $\alpha < 1$, it is decreasing in the ball $B(0;1/3)$.  
Hence $\omega^{\alpha}$ has a positive lower bound, say $c$. Next, using the notation from Section 4, 
we estimate $\int A_{2,\kappa}(\omega^{\alpha}_+ ) \,dx$. Since $A_{2,\kappa}$ is nondecreasing we have 
\[ 
\int A_{2,\kappa}(\omega^{\alpha}_+ ) \,dx \leq \int |\omega^{\alpha}|^2 \log^{2\kappa}(\omega^{\alpha} + 2)\, dx   
\] 
\[\leq C(\|\omega^{\alpha}\|_{L^2}) \int_{B(0;1/3)} \frac{1}{|x|^2|\log|x||^{2\alpha}} 
\left| \log^{2\kappa} \frac{1}{|x||\log|x||^{\alpha}} \right| \, dx,\] 
using the fact that $|\omega^{\alpha}| \geq c>0$ on $B(0;1/3)$, 
\[ \leq C \int_{B(0;1/3)} \frac{1}{|x|^2 |\log|x||^{2\alpha}}|\log|x||^{2\kappa} \, dx = 
C\int_0^{1/3} \frac{1}{r|\log r|^{2\alpha - 2\kappa}} \, dr < \infty, 
\] 
as long as $2\alpha - 2\kappa > 1$, i.e., $\kappa < \alpha - 1/2$. The last inequality is due to the fact that the double logarithm grows slower than the single logarithm.  
 
The condition that $\kappa \geq 0$ arises from the definition of the Zygmund spaces.  

\end{proof}

\begin{theorem} \label{limitcase} 
If $\alpha <1$ then $u^{\alpha}_+ \equiv K \ast \omega^{\alpha}_+$ is 
unbounded. If $1/2<\alpha \leq 2/3$ then 
$u^{\alpha}_+ |\omega^{\alpha}_+|^2$ is not locally integrable.
\end{theorem}

\begin{proof}
   
        We will show that the first component of $u^{\alpha}_+$, which we denote by $u_1$, is greater 
than or equal to $C |\log|x||^{1-\alpha}$ in a suitably small
neighborhood of the origin. It is easy to see that this result 
proves both assertions in the statement of the theorem. 
 
        First we compute $u_1$ on the horizontal axis. Note that $\omega^{\alpha}_+$ is even with respect to 
$x_1=0$. Then $u_1$ has the same symmetry, due to the specific form of the Biot-Savart kernel,  and in particular $u_1(x_1,0) = u_1(-x_1,0)$. Therefore, 
it is enough to compute $u_1(x_1,0)$ for $x_1>0$. We have 
\[u_1(x_1,0) = \int_{B^+(0;1/3)} \frac{y_2}{2\pi|x-y|^2}\,\frac{1}{|y||\log|y||^{\alpha}}\,dy \] 
 
\[ = \frac{1}{2\pi}\int_0^{1/3} \int_0^{\pi}\frac{r\sin \theta}{(x_1-r\cos\theta)^2+ (r\sin \theta)^2}\,d\theta 
\frac{1}{|\log r|^{\alpha}}\,dr, \] 
after changing to polar coordinates. Explicitly evaluating the integral in $\theta$ and subsequently  
implementing the change of variables $s=r/x_1$ we find 
\[2 \pi u_1(x_1,0) = \int_0^{1/3} \frac{1}{|\log r|^{\alpha}}\,\frac{1}{x_1}\log\left|\frac{r+x_1}{r-x_1}\right|\,dr\] 
 
\[= \int_0^{1/(3x_1)} \frac{1}{|\log s x_1|^{\alpha}} \log\left| \frac{s+1}{s-1} \right| \, ds\] 
 
\[= \int_0^2 \frac{1}{|\log s x_1|^{\alpha}} \log\left| \frac{s+1}{s-1} \right| \, ds + \int_2^{1/(3x_1)} 
\frac{1}{|\log s x_1|^{\alpha}} \log\left| \frac{s+1}{s-1} \right| \, ds \]
\[\equiv \mathcal{I} + \mathcal{J}.\] 
 
 We assume $0 \leq x_1 < 1/6$ and we estimate $\mathcal{I}$: 
 
\[0 \leq \mathcal{I} \leq \frac{1}{(\log 3)^{\alpha}} \int_0^2 \log \left| \frac{s+1}{s-1} \right| \, ds \equiv C <\infty.\] 
 
Next we estimate $\mathcal{J}$ from below. We begin with two observations. For $2 < s < 1/(3x_1)$ we have: 
 
\[  \frac{1}{|\log sx_1|^{\alpha}} \geq  \frac{1}{|\log 2x_1|^{\alpha}}; \] 
and  
\[ \log \frac{s+1}{s-1} > \frac{1}{s}. \] 
 
Therefore,  
\[\mathcal{J} \geq \frac{|\log 6x_1|}{|\log{2x_1}|^{\alpha}} \geq \frac{1}{2} |\log x_1|^{1-\alpha}, \] 
where the last inequality was derived assuming further that $x_1 \leq 1/36$. 
 
In summary, we have shown that  
\begin{equation} \label{nearor} u_1(x_1,0) \geq C|\log|x_1||^{1-\alpha} \;\;\text{ if } |x_1| \leq \frac{1}{36}, 
\end{equation} 
for some $C>0$. In addition, it follows from the specific form of the Biot-Savart law that $u_1(x_1,0) \geq 0$ for all $x_1$. 
 
Recall the radially symmetric function $\omega^{\alpha}$, introduced in \eqref{cizania}.  Consider the vorticity  
$\omega^{\alpha} - \omega^{\alpha}_+$, supported in the lower half-plane. Let  
\[v_1=v_1(x_1,x_2) = \int \frac{y_2 - x_2}{2\pi|x-y|^2}(\omega^{\alpha} - \omega^{\alpha}_+)(y)\,dy,\] 
be the first component of the associated velocity. Then $v_1$ is a harmonic function in the upper half-plane, whose boundary value, by symmetry, is equal to $-u_1(x_1,0)$, since the horizontal velocity associated to 
$\omega^{\alpha}$ vanishes on the horizontal axis. We may thus write, using the Poisson kernel for the upper half-plane, 
\begin{equation} \label{closure} 
v_1(x_1,x_2) = -\frac{1}{\pi}\int_{-\infty}^{\infty} \frac{x_2\, u_1(s,0)}{(x_1-s)^2 + x_2^2} \, ds, \, \text{ if }  
x_2>0. 
\end{equation} 
 
Note that $u^{\alpha} \equiv v_1 + u_1$ is the velocity associated to $\omega^{\alpha}$. In view of Lemma \ref{symbd} we have that $u^{\alpha}$ is bounded and there exists $C>0$ such that: 
\[\|u^{\alpha}\|_{L^{\infty}} \leq C\|\omega^{\alpha}\|_{L^2}.\]  
 
In what follows we will show that $v_1 \leq -C|\log |x||^{1-{\alpha}}$ for sufficiently small $|x|$, with $x_2>0$. By virtue of the previous observation this is enough to conclude the proof.  
 
Let $0< \delta < 1/36$. Using \eqref{nearor} and the fact that $u_1$
is nonnegative on $x_2 = 0$, we find  for $x_2>0$, 
\[v_1(x_1,x_2) \leq -\frac{1}{\pi}\int_{-\delta}^{\delta} \frac{x_2\, u_1(s,0)}{(x_1-s)^2 + x_2^2} \, ds \] 
\[\leq -C|\log \delta|^{1- \alpha} \left\{ \arctan \left( \frac{x_1 + \delta}{x_2} \right) -  
\arctan \left( \frac{x_1 - \delta}{x_2} \right)\right\},\] 
by explicitly integrating the Poisson kernel in the interval $(-\delta,\delta)$. 
 
Next, let $x = (\delta/2)(\cos \theta, \sin \theta)$ with $0\leq \theta \leq \pi$. Then  
\[v_1(x) \leq -C |\log 2|x||^{1-\alpha} \left[ \arctan \left( \frac{\cos \theta + 2}{\sin \theta} \right) -  
\arctan \left( \frac{\cos \theta - 2}{\sin \theta} \right)\right]\] 
\[ \equiv -C |\log 2|x||^{1-\alpha}g(\theta).\] 
It is easy to compute the minimum of $g(\theta)$, thereby verifying that $g(\theta) \geq 2\arctan 2 > 0$  
for all $\theta \in [0,\pi]$. We have therefore shown that, 
for any $x=(x_1,x_2)$ with $x_2 > 0$ and $|x| \leq 1/72$,  $v_1(x_1,x_2) \leq -C |\log 2|x||^{1-\alpha}$. The conclusion follows as $|\log 2|x||\geq (1/2)|\log |x||$ for any $|x| < 1/4$. 
\end{proof} 
 
\begin{remark} We emphasize that we have proved above that there exist constants $C>0$, $0<r_0<1/72$ such that 
\begin{equation} \label{saco} 
u_1(x) \geq C |\log |x||^{1-\alpha}, \,\mbox{ for } x \in B^+(0;r_0). 
\end{equation} 
\end{remark} 
 
\vspace{0.5cm} 
 
We wish to use Lemma \ref{precreg} and Theorem \ref{limitcase} to  
draw two separate conclusions. The first is that $L^2(\log L)^{1/4}$ contains vorticities   
whose associated velocities are unbounded. Indeed, it is enough to consider $\omega^{\alpha}_+$,  
for $3/4<\alpha<1$. The second conclusion is that there are difficulties in making sense, as a distribution, of 
$u |\omega|^2$ for an arbitrary vorticity in $L^2$. In fact, we have already shown that $u^{\alpha}_+ 
|\omega^{\alpha}_+|^2$ is not locally integrable if $1/2 < \alpha \leq 2/3$.  
From Lemma \ref{precreg} it follows that  
$\omega^{\alpha}_+ \in L^2(\log L)^{\kappa}$, for some $0\leq\kappa<1/6$ if $1/2 < \alpha \leq 2/3$.  Although suggestive, the non-integrability of $u^{\alpha}_+ |\omega^{\alpha}_+|^2$ does not exclude the possibility that  
$u^{\alpha}_+ |\omega^{\alpha}_+|^2$ gives rise to a well-defined distribution. 
One may recall the way in which the non-integrable functions $1/s$ and $1/s^2$ can be  
identified with the distributions pv-$1/s$ and pf-$1/s^2$. 
 
We must address more closely the problem of identifying $u|\omega|^2$ with a distribution.  
In view of Definition \ref{uoomega} one might
suspect that by re-arranging the Biot-Savart law in a clever way and
using the antisymmetry of the kernel, it would be possible to give meaning
to $u|\omega|^2$ in a consistent manner, even if $\omega$ is only in $L^2_c$.    
The antisymmetry of the Biot-Savart kernel  
has been used on more than one occasion to prove results of this nature; for instance it  
was used to define the nonlinear term $u \cdot \nabla \omega$, when  
$\omega \in {\mathcal BM}_+ \cap H^{-1}_{\loc}$, by S. Schochet in \cite{Sch}.  
We will see that this strategy would not be successful in this case. 
 
Ultimately, our purpose here is to examine the sharpness of the condition $\omega \in L^2(\log L)^{1/4}$, which 
we showed to be sufficient to define the term $u|\omega|^2$. This condition  
was used in Theorem \ref{ensbalance} and Definition \ref{uoomega}. We would like to 
argue through a counterexample that it is not possible to make sense of $u |\omega|^2$ for 
arbitrary $\omega \in L^2(\log L)^{\kappa}$, with $0\leq\kappa < 1/6$.    If we wish to attribute meaning to  
$u |\omega|^2$ (as a distribution) for any $\omega \in X \subseteq L^2$, then the key issue is the nature of the nonlinear map $T:\omega \mapsto u|\omega|^2$, from $X$ to $\mathcal{D}^{\prime}$.  
First, note that $T$ is well-defined for $X=L^p_c$, $p>2$, since then $u = K\ast \omega$ is bounded.  
Next, note that Definition \ref{uoomega} actually consists of the continuous extension of $T$ to  
$X=(L^2(\log L)^{1/4})_c$. We will show through the counterexample we present  that there  
is no continuous extension of $T$ from $L^p_c$, $p>2$ to $X=(L^2(\log L)^{\kappa})_c$, $0\leq\kappa<1/6$, and hence, to $X=L^2_c$. In fact we will prove that our example $\omega^{\alpha}_+$, with $1/2<\alpha \leq 2/3$, can
be approximated in $(L^2(\log L)^{\kappa})_c$, $0\leq\kappa<\alpha-1/2$, by a sequence $\omega^n_+ \in L^{\infty}_c$
for which $\int u^n_+|\omega^n_+|^2 \to \infty$ as $n \to \infty$, thereby reaching the desired conclusion.

\begin{theorem} \label{finitecase}
Let $x = (x_1,x_2)$ with $x_2 \geq 0$. Fix $1/2 < \alpha \leq 2/3$. For each $n \in \N$ we define the approximate vorticity by: 
\begin{equation} 
\omega^n_+ = \omega^n_+(x)= \left\{ 
\begin{array}{l} 
\omega^{\alpha}_+ (x) \, \text{ if } |x| > 1/n,\\ 
\\ 
\frac{n}{|\log n|^{\alpha}} \, \text{ if } |x| \leq 1/n. 
\end{array} \right. 
\end{equation} 
Then $\omega^n_+ \to \omega^{\alpha}_+$, as $n \to \infty$, strongly in $L^2(\log L)^{\kappa}$ for all $0\leq\kappa < \alpha-1/2$.  
 
Denote $u^n_1$ the first component of $K\ast\omega^n_+$. Then it also holds that  
\begin{equation} \label{yeah} \lim_{n\to +\infty}\int u^n_1|\omega^n_+|^2 \, dx = +\infty. \end{equation}
\end{theorem}

\begin{proof} 
Our first step is to show that $u^n_1$ is nonnegative in $B^+(0;r_0)$, if $n$ is large enough, 
where $r_0$ is such that \eqref{saco} holds. We require two different arguments, one for  
$|x| \leq 2/n$ and another for $2/n<|x|<r_0$. We will begin with the latter. 
 
Let $W_n = \omega^{\alpha}_+ - \omega^n_+ \geq 0$, which is a function with support in $B^+(0;1/n)$. 
Let $e_n$ be the first component of $K \ast W_n$, i.e., the error in the velocity induced by the truncation. 
Therefore, $u^n_1 = u_1 - e_n$.  It follows from \eqref{saco} that  
\begin{equation} \label{oestt} u^n_1(x) \geq C|\log|x||^{1-\alpha} - e_n(x), \, \mbox{ for } x \in B^+(0;r_0).  
\end{equation} 
We will prove that  
\begin{equation} \label{estt} |e_n(x)| \leq C/(\log n)^{\alpha}, \, \mbox{ for } |x| > 2/n. \end{equation} 
For $x \in B^+(0;r_0)$, $|x|>2/n$ we estimate: 
\begin{equation} \label{thing1} 
|e_n(x)| \leq \int_{B^+(0;1/n)} \frac{1}{|x-y|} W_n(y) \,dy \leq Cn \int_{B^+(0;1/n)} W_n(y) \,dy, \end{equation} 
as $|x-y| \geq 1/n$, 
\begin{equation} \label{thing2} 
 = Cn \int_0^{1/n} \left( \frac{1}{|\log r|^{\alpha}} 
- \frac{nr}{(\log n)^{\alpha}} \right) \, dr, \end{equation} 
after changing to polar coordinates,  yielding \eqref{estt}.  
 
As $|\log|x||^{1-\alpha}$ is decreasing with respect to $|x|$, it follows from \eqref{oestt} and \eqref{estt} 
that one can choose $n_0$ sufficiently large so that if $n>n_0$ and $|x| > 2/n$, with $x \in B^+(0;r_0)$, then  
$u_1^n(x) \geq 0$.

Now we address the case $|x| \leq 2/n$. We will show that  
\begin{equation} \label{creative} 
u_1^n(x) \geq C (\log n)^{1-\alpha} 
\end{equation}  
for $x$ in this region. 
The proof closely parallels the proof of Theorem \ref{limitcase}. We begin by estimating $u_1^n(x_1,0)$ if $|x_1|< 2/n$. We have: 
\[\pi u_1^n(x_1,0) = \int_{1/n}^{1/3} \frac{1}{r|\log r|^{\alpha}} \frac{r}{|x_1|} \log \left|   
\frac{(r/|x_1|)+1}{(r/|x_1|)-1}\right| \, dr\]  
\[ + \frac{n}{(\log n)^{\alpha}} \int_0^{1/n}\frac{r}{|x_1|} \log \left|   
\frac{(r/|x_1|)+1}{(r/|x_1|)-1}\right| \, dr \] 
\[\geq \int_{2/n}^{1/3} \frac{1}{r|\log r|^{\alpha}} g(r/|x_1|) \, dr,\] 
where $g(s) \equiv s\log\left|(s+1)/(s-1)\right|$. It can be easily verified that $g(s) > 1$ if $s>1$. Therefore, as $r/|x_1| > 1$ for $r>2/n$ and $|x_1|<2/n$, we obtain 
\begin{equation} \label{ifyouwant} 
\pi u_1^n(x_1,0) \geq \int_{2/n}^{1/3} \frac{1}{r|\log r|^{\alpha}}\, dr \geq C|\log n|^{1-\alpha}, 
\end{equation} 
for $n$ sufficiently large. We also know that $u^n_1(x_1,0) \geq 0$ for all $x_1$.  
 
Let $\omega^n$ be the radially symmetric extension of $\omega^n_+$ and set $v^n_1$ to be the first component of  
$K\ast (\omega^n - \omega^n_+)$. As in the proof of Theorem \ref{limitcase}, we find that 
\[\pi v^n_1(x_1,x_2) = - \int_{-\infty}^{+\infty} \frac{x_2 \, u^n_1(s,0)}{x_1 - s)^2 + x^2_2} \, ds \] 
\[ \leq - \int_{-2/n}^{2/n} \frac{x_2 \, u^n_1(s,0)}{x_1 - s)^2 + x^2_2} \, ds \] 
\[\leq - C(\log n)^{1-\alpha} \left[ \arctan\left( \frac{x_1 + (2/n)}{x_2} \right) -  
\arctan\left( \frac{x_1 - (2/n)}{x_2} \right)   \right] ,\] 
by \eqref{ifyouwant}. It is easy to see that, if $|x_1| < 2/n$ and $0< x_2 < 2/n$, then the difference of arctangents above is bounded from below by $\arctan 1 = \pi/4$. Therefore we deduce that, if $n$ is sufficiently large, then  
$v^n_1(x) \leq - C(\log n)^{1-\alpha}$ for $x \in B^+(0;2/n)$. Then, as in the proof of Theorem \ref{limitcase}, we obtain \eqref{creative} as long as $n$ is large enough. This completes the proof that $u_1^n$ is 
nonnegative in $B^+(0;r_0)$ for $n$ large enough.  
 
Let $\mathcal{U}_n \equiv B^+(0;r_0) \setminus B^+(0;1/\sqrt[3]{n})$. Recall that $e_n$ is the error  
in the first component of velocity, due to truncation. We will show that there exists $C>0$, such  
that for $n$ sufficiently large we have 
\begin{equation} \label{stpost} 
\left|\int_{\mathcal{U}_n} e_n |\omega^n_+|^2 \, dx \right| \leq  \frac{C}{(\log n)^{\alpha}}. 
\end{equation}   
 
In fact, we observe first that for $x \in \mathcal{U}_n$ we have 
\[|e_n(x)| \leq C \sqrt[3]{n} \int_{B^+(0;1/n)} W_n(y) \, dy, \]  
as $|x-y| \geq 1/(2\sqrt[3]{n})$ for $n$ sufficiently large and $|y| \leq 1/n$, so that
\[|e_n(x)|\leq \frac{C}{n^{2/3} (\log n)^{\alpha}}, \] 
as in the proof of \eqref{estt}, see \eqref{thing1}, \eqref{thing2}. 
Additionally, for $x \in \mathcal{U}_n$, 
\[|\omega_+^n(x)|^2 = \frac{1}{|x|^2|\log |x||^{2\alpha}} \leq C(r_0) n^{2/3}.\] 
Estimate \eqref{stpost} follows immediately from these two observations. 
 
We now complete the proof of \eqref{yeah}. 
We note that
\[\int_{B^+(0;r_0)} u^n_1|\omega^n_+|^2 \,dx \geq \int_{\mathcal{U}_n} u^n_1|\omega^n_+|^2 \,dx ,\] 
as $u^n_1 \geq 0$ in $B^+(0;r_0)$, 
\[= \int_{\mathcal{U}_n} u_1|\omega^{\alpha}_+|^2 \,dx - \int_{\mathcal{U}_n} e_n |\omega^n_+|^2 \,dx \equiv  
\mathcal{I}_n + \mathcal{E}_n,\] 
where we have used that $\omega^n_+ = \omega^{\alpha}_+$ in $\mathcal{U}_n$ and $u^n_1 = u_1 - e_n$. 
 
By \eqref{stpost} we obtain that $\mathcal{E}_n \to 0$ as $n \to \infty$. 
Moreover, we have established in Theorem \ref{limitcase} that  
\begin{equation} \label{dirty} 
\int_{B^+(0;r_0)} u_1(x)|\omega^{\alpha}_+(x)|^2 \, dx = \infty. 
\end{equation} 
Therefore, using the Monotone Convergence Theorem, we find that $\mathcal{I}_n \to \infty$ as $n \to \infty$.  
We conclude that  
\[\lim_{n \to \infty} \int_{B^+(0;r_0)} u^n_1|\omega^n_+|^2 \,dx = \infty.\]  
To finish the proof of \eqref{yeah} we observe that arguments similar to those used above imply that $u^n_1$ is bounded in $B^+(0;1/3)\setminus B^+(0;r_0)$; the same is true of $\omega^n_+$ by construction. This completes the proof of \eqref{yeah}.  
 
Finally, we turn to the convergence of $\omega^n_+$ to $\omega^{\alpha}_+$. Let $0 \leq \kappa < \alpha - 1/2$.  
We estimate the difference in the Zygmund space $L^2(\log L)^{\kappa}$. We have that: 
\begin{equation}\label{Zygest1} 
 \int A_{2,\kappa}\left(\frac{W_n}{\|W_n\|_{2,\kappa}} \right)\, dx = 1,\end{equation} 
since $W_n$ does not vanish identically.  
We observe that $0 \leq W_n \leq \omega^{\alpha}_+ \chi_{B^+(0;1/n)}$. By Lemma \ref{precreg},   
$\omega^{\alpha}_+ \in L^2(\log L)^{\kappa}$. Therefore, since $A_{2,\kappa}$ is nondecreasing, it follows that 
\begin{equation} \label{Zygest2} 
\int A_{2,\kappa}(W_n)\, dx \leq \int_{B^+(0;1/n)} A_{2,\kappa}(\omega^{\alpha}_+ )\, dx \to 0,\end{equation} 
as $n \to \infty$ by continuity of integrals. Now, recall that $A_{2,\kappa}$ is convex. Therefore, 
\begin{equation} \label{Zygest3} 
A_{2,\kappa}\left(\frac{W_n}{\|W_n\|_{2,\kappa}} \right) \leq \frac{1}{\|W_n\|_{2,\kappa}} A_{2,\kappa}(W_n). 
\end{equation} 
By virtue of \eqref{Zygest1} and \eqref{Zygest3} we find  
\[\|W_n\|_{2,\kappa} \leq \int A_{2,\kappa}(W_n)\, dx.\] 
Using \eqref{Zygest2} then  implies that $\|W_n\|_{2,\kappa} \to 0$ as we wished. 

\end{proof}
 
    We emphasize at this point that this section was concerned with the cubic nonlinearity $u|\omega|^2$ without  
reference to dynamics. Something strange might occur with enstrophy
dissipation and with the transport enstrophy defect at the initial
time for a weak solution of incompressible 2D Euler obtained with
$\omega^{\alpha}_+$ as initial data. We do not offer any prognosis, as the answer depends on how the initial snarl in the term $u|\omega|^2$ would resolve itself for positive time. It would be very interesting to determine what happens, but this problem 
does not seem tractable. 
 
\section{Counterexample for Eyink's conjecture}   
 
In this section we will present a counterexample to Eyink's conjecture, as formulated in Section  
3. We will exhibit a family of solutions to the 2D Navier-Stokes
equations, which converge, as viscosity vanishes, to a stationary solution of 2D Euler. This stationary solution is such that both $Z^T$ and $Z^V$ exist and $Z^T$ vanishes identically while $Z^V$ does not.

We consider $\omega_{0}$ of the form:  
\begin{equation} \label{e:data}  
  \omega_{0} (x)=\phi(x)\frac{1}{|x|}, \quad \phi\in C^{\infty}_c
(\real^{2}),  
\end{equation}  
with $\phi$ radially symmetric, $\supp\phi \subset B(0;1)$, $\phi\equiv 1$ on   
$B(0;1/2)$. Note that such $\omega_{0}$ belongs to $L^{2,\infty} \cap
L^{p}$, $1\leq p<2$.

It is well known that any radially symmetric vorticity configuration  
$\omega = \omega(x) = \rho(|x|)$ gives   
rise to an exact steady solution $u$ of the incompressible Euler equations, see \cite{majbert}. As in 
Lemma \ref{symbd}, the 2D Biot-Savart law becomes:  
\begin{equation} \label{e:radial}  
       u(x) = \frac{x^{\perp}}{|x|^2}\int_{0}^{|x|}  s\rho(s)\,ds.  
\end{equation}  
Such steady solutions are called Rankine vortices.     
  
\begin{remark} If $\phi$ is chosen instead so that $\,\int \omega_0(x)  
\,dx=0$, then $u$ defined in  \eqref{e:radial} is  compactly  supported,  
vanishing outside $\supp\phi$ (see \cite{DiPM1}).  This observation would 
allow us to adapt the present example to the periodic case.  
\end{remark}  
  
Similarly, if $\omega_{\nu}$ is the solution of the {\em heat} equation  
\begin{equation} \label{heet}  
\partial_t \omega_{\nu} = \nu \Delta \omega_{\nu}, \end{equation} 
with radially symmetric initial data $\omega_0$, then $u_{\nu} \equiv K \ast \omega_{\nu}$ 
is a solution  of the 2D Navier-Stokes equations with initial
vorticity $\omega_0$ and viscosity $\nu$. 
  
We will show that $\omega_{0}$ belongs to $B^{0}_{2,\infty}$
and that the sequence $\omega_{\nu}$ satisfies the hypothesis of the
Eyink's Conjecture, however we postpone the proof of this fact to the end of this section, 
see Proposition \ref{Besov}.
 
In what follows, we recall the notation used in Section 3. If $j_{\epsilon}$ is a (radially 
symmetric) Friedrichs mollifier, then we denote  $j_{\epsilon} \ast
\omega_0$ with $\omega_{\epsilon}$.   
We introduce the approximate transport defect $Z_{\epsilon}(\omega_0)$ and the approximate viscous defect  
$Z^{\nu}(\omega_{\nu})$ as defined in Section 3. 

We state below the main result of this section. 
  
\begin{theorem} \label{t:example}  
The enstrophy defects $Z^T(\omega_0)$ and $Z^V(\omega_0)$ both exist. Moreover,  
\[ Z^T(\omega_0) \equiv 0 \,\, \mbox{ while } \, Z^V(\omega_0) = \frac{4\pi^3}{t} \delta_0, \] 
where $\delta_0$ is the Dirac measure supported at the origin.  
\end{theorem}  
 
\begin{proof} 
  
To prove that $Z^T(\omega_0)$ exists and vanishes identically we observe that  
$\omega_{\epsilon}$ remains radially symmetric by construction and the flow  
lines of $u_{\epsilon} = j_{\epsilon} \ast K \ast \omega_0$ are concentric circles  
centered at the origin. Therefore we find  
\[Z_{\epsilon}(\omega_0) = -\nabla\omega_{\epsilon}\cdot  \left(   
  (u\omega_{\epsilon})_{\epsilon}-u_{\epsilon} \omega_{\epsilon}\right)= 0,\] 
so that $Z^T(\omega_0) \equiv 0$. 
  
In the rest of the proof, we will discuss the viscous enstrophy defect. We begin   
by deriving sharp asymptotic estimates for   
$\nu\|\nabla \omega_{\nu}\|_{L^{2}}^2=\|Z^{\nu}(\omega_{\nu})\|_{L^{1}}$.    
This is accomplished in the following proposition. 
  
\begin{proposition} \label{p:normest}  
For each $t>0$, the approximate viscous enstrophy defect satisfies:
\begin{equation} \label{e:normest}  
    \|Z^{\nu}(\omega_{\nu})\|_{L^{1}} = \frac{4\pi^3}{t} + \text{o}(1),  
\end{equation}  
as $\nu\to 0^{+}$.  
 \end{proposition}  
     
\begin{proof}[Proof of Proposition.]  
By  Plancherel's Theorem we have 
\begin{equation} \label{e:normest1}  
  \|Z^{\nu}(\omega_{\nu})\|_{L^{1}}=\nu\int_{\real^{2}}|\xi|^{2}\,e^{-t\nu|\xi|^{2}}  
   \,|\Hat{\omega}_{0}(\xi)|^2\,d\xi.  
\end{equation}   

We begin by estimating the Fourier transform   
of $\omega_{0}$. Set  
\[e=e(\xi) = |\xi||\Hat{\omega}_0(\xi)| - 2\pi.\] 
We will show that $e$ is a bounded function which vanishes along rays near $\infty$, i.e., for each $\xi \neq 0$ fixed, $|e(s\xi)| \to 0$ as $s \to \infty$. To this end, fix $\xi \neq 0$ and write $\xi = r \sigma$, with  
$|\sigma| = 1$ and $r = |\xi|$. 
We recall that
\[\left( \frac{1}{|x|} \right)^{\Hat{}}(\xi) = \frac{2\pi}{|\xi|},\] 
(see Lemma 1 of Chapter V of \cite{St70} for details)
and hence 
\begin{equation} \label{hatomega0} 
\Hat{\omega}_0 (\xi)= \left(\frac{2\pi}{|z|} \ast \Check{\phi}(z) \right)(-\xi),\end{equation} 
by the usual properties of the Fourier transform. As $\phi \in C^{\infty}_c$ it follows that $\Check{\phi} \in \Sc$, the Schwartz space of rapidly decaying smooth functions.  
 
Using \eqref{hatomega0} now gives: 
\[e(s\xi) =|s||\xi||\Hat{\omega}_0(s\xi)| - 2\pi = 
|sr| \left| \left(\frac{2\pi}{|z|} \ast \Check{\phi}(z) \right)(-s\xi) \right| - 2\pi\] 
\[= |sr| \left| \int \frac{2\pi}{|y|}\Check{\phi}(-s\xi - y) \, dy 
\right| - 2\pi\] 
\[=|sr| \left| \int \frac{2\pi}{|sr||z|}\,\Check{\phi}(-|sr|\sigma - |sr|z) \, |sr|^2dz 
\right| -2\pi,\] 
after making the change of variables $y = |sr| z$, 
\[=\left| \int \frac{2\pi}{|z|}\,|sr|^2\Check{\phi}(|sr|(-\sigma - z)) \, dz 
\right| -2\pi \] 
\[\equiv \left| \int \frac{2\pi}{|z|}\Check{\phi}_{|sr|}(-\sigma - z) \, dz 
\right| -2\pi =   
\left|\left(\frac{2\pi}{|z|}  \ast \Check{\phi}_{|sr|}(z) \right)(-\sigma)\right| - 2\pi,\] 
where $\Check{\phi}_M (\cdot)\equiv M^2 \Check{\phi} (M\cdot)$.  
 
It is easy to see that $|e(s\xi)|$ is uniformly bounded in both $s$
and $\xi$, since $C/|z|$ is a locally integrable function, bounded
near infinity, and $\Check{\phi}_{|sr|}$ is small near infinity and
integrable with constant integral with respect to $sr$. We simply
estimate the convolution above by distinguishing points $z$ near $\sigma$
and points  $z$ far from  
$\sigma$. Then we use the fact that $|\sigma| = 1$. 
 
Note that $\Check{\phi}_{|sr|} \to \delta_0$ in $\mathcal{D}^{\prime}$
as $s \to \infty$ so that the convolution above should, in principle,
converge to $2\pi/|-\sigma | = 2\pi$. The difficulty in making this
argument precise is that 
neither $2\pi / |z|$ nor $\Check{\phi}$ are compactly supported. 
 
Let $0<\beta < 1/4$. Set $\eta_{\beta} = \eta_{\beta}( \xi)$ a radially symmetric smooth cut-off function of the ball of radius $\beta$, so that $\eta_{\beta}$ is identically $1$ in $B(0;\beta)$ and vanishes in  
$\R^2 \setminus B(0;2\beta)$. We use $\eta_{\beta}$ to write: 
\[\left(\frac{2\pi}{|z|}  \ast \Check{\phi}_{|sr|}(z) \right)(-\sigma) = \] 
\[= \int \frac{2\pi}{|-\sigma-z|} \eta_{\beta}(z)\Check{\phi}_{|sr|}(z) \, dz 
+\int \frac{2\pi}{|-\sigma -z|}(1- \eta_{\beta}(z))\Check{\phi}_{|sr|}(z) \, dz \] 
\[\equiv \mathcal{J}_1 + \mathcal{J}_2.\] 
 
Note that, for each fixed $\sigma$, with $|\sigma|=1$, the function  
$2\pi\eta_{\beta}(z) /|-\sigma-z| $ is smooth and compactly supported, which implies that $\mathcal{J}_1 \to 2\pi$ as $s \to \infty$. We show that $\mathcal{J}_2 \to 0$: 
\[\left| \int \frac{2\pi}{|-\sigma -z|}(1- \eta_{\beta}(z))\Check{\phi}_{|sr|}(z) \, dz \right| \leq \] 
\[\leq \left| \int_{\beta<|z|<2} \frac{2\pi}{|-\sigma -z|} 
\Check{\phi}_{|sr|}(z) \, dz \right|+ \left| \int_{|z|>2}  
\frac{2\pi}{|-\sigma -z|}\Check{\phi}_{|sr|}(z) \, dz \right|\] 
\[\leq   \|\Check{\phi}_{|sr|}\|_{L^{\infty}(\{|z|>2\beta\})}\left| \int_{2<|z|<2\beta} \frac{2\pi}{|-\sigma -z|} \, dz  \right| 
+ 2\pi \left| \int_{|z|>2}\Check{\phi}_{|sr|}(z) \, dz \right|.\] 
Clearly each term above vanishes as $s \to \infty$. 
 
Finally, we may now write: 
\begin{equation} \label{e:square} 
 |\xi|^2|\Hat{\omega}_0(\xi)|^2 = 4\pi^2 + 4\pi e(\xi) + |e(\xi)|^2, 
\end{equation}
so that, from \eqref{e:normest1}, we find 
\[\|Z^{\nu}(\omega_{\nu})\|_{L^{1}}=\nu\left( 
\int_{\real^{2}}4\pi^2\,e^{-t\nu|\xi|^{2}}\,d\xi + \int_{\real^{2}}(4\pi e(\xi) + |e(\xi)|^2)\, 
e^{-t\nu|\xi|^{2}}\,d\xi \right)\] 
\[=\frac{1}{t}\int_{\real^{2}}4\pi^2\,e^{-|z|^2}\,dz + \frac{1}{t}\int_{\real^{2}}\left( 4\pi e \left(\frac{z}{\sqrt{t\nu }}\right) + \left|e\left(\frac{z}{\sqrt{t\nu }}\right)\right|^2 \right)\,e^{-|z|^2}\,dz .\] 
Since we have already shown that $e=e(z)$ is a bounded function and
that $\lim_{\nu \to 0^+} e(z/\sqrt{t\nu}) = 0$, we deduce
using the Dominated Convergence Theorem that  
\[\lim_{\nu\to 0^+}\|Z^{\nu}(\omega_{\nu})\|_{L^{1}} = \frac{4\pi^3}{t},\] 
as we wished. 
 
\end{proof}  
  
In view of the proposition above we find that, for each fixed $t>0$, the set  
$\{Z^{\nu}(\omega_{\nu}),\;\nu>0\}$  
is uniformly bounded in $L^{1}$. Therefore,  using the Banach-Alaoglu Theorem, for each $t>0$ there is  
a sequence converging weakly to a Radon measure.  Each of these measures is, in fact,  
a multiple of the Dirac measure, $C(t)\,\delta_0$, by virtue of the following   
claim, which we will prove later.    
  
\begin{claim}  
  Any converging sequence of  $\{Z^{\nu}(\omega_{\nu}),\;\nu>0\}$  
  converges to a  distribution supported at the origin.  
\end{claim}  
  
Given the Claim we may conclude that $Z^{\nu}(\omega_{\nu})$ is itself convergent (to a  positive   
measure). To establish this result, it is enough to show that  $C(t)$ is independent of   
the particular sequence $Z^{\nu_{k}}(\omega_{\nu_k})$. To this end, we fix a converging subsequence  
$Z^{\nu_k}(\omega_{\nu_k})$. We begin by observing that $Z^{\nu}(\omega_{\nu})$ is a {\it tight} family  
of functions in $L^1$ with respect to the parameter $\nu$. Indeed,
$Z^{\nu}(\omega_{\nu})= \nu |\nabla \omega_{\nu}|^2$ and
$\omega_{\nu}$ is the convolution of a compactly supported function
with the heat kernel, so it is immediate to verify that $\int_{|z|>M}
Z^{\nu}(\omega_{\nu}) \, dz \to 0$ as $M \to \infty$, uniformly in
$\nu$. Fix now $\epsilon > 0$ and choose $M$ so large that $0< \int_{|z|>M} Z^{\nu}(\omega_{\nu}) \, dz < \epsilon$. Then, if $\psi_M$ is a smooth cut-off of the ball of radius $M+1$, we have 
\[ \left| \int_{\real^{2}} Z^{\nu_{k}}(\omega_{\nu_k})\,dx - \int_{\real^2}   
   Z^{\nu_{k}}(\omega_{\nu_k}) \,\psi_M\,dx\right|= \] \[=\left|\int_{\real^2}   
   Z^{\nu_{k}}(\omega_{\nu_k}) \,(1-\psi_M)\,dx\right| < \epsilon,  \]  
By Proposition \ref{p:normest} and the Claim the left-hand side  
converges to $4\pi^3 / t - C(t)$ as $k \to \infty$. As $\epsilon$ is arbitrary it follows that  
$C(t) = 4\pi^3/t$, independent of the sequence $\nu_k$, as desired.
  
In summary, we have deduced that  
\begin{equation}  
    \lim_{\nu\to0^{+}} Z^{\nu}(\omega_{\nu}) =\frac{4\pi^3}{t} \delta_0, \text{ in } \mathcal{D}^{\prime}.  
\end{equation}  
It remains to establish the Claim.  
   
\begin{proof}[Proof of Claim]  
 We prove that, for any $\eta>0$ and any $f \in C^{\infty}_c$ with $\supp f \subset \real^{2}\setminus B(0;\eta)$, we 
have 
\begin{equation} \label{e:claim1}  
     \lim_{\nu\to0^{+}} \int _{\real^{2}} Z^{\nu}(\omega_{\nu})\, f\,dx=0.  
\end{equation}  
 
The proof involves a simple estimate on $\omega_{\nu}$. 
Let $H_{\nu}=H_{\nu}(x,t)=(4\pi \nu t)^{-1}\,e^{-|x|^{2}/(4\nu t)}$
denote the heat kernel in $\real^{2}$.
Recall that $\omega_{\nu}$ satisfies \eqref{heet} so that we 
may write $\omega_{\nu} = H_{\nu} \ast \omega_0$. Fix $\eta > 0$ and let $\varphi \in C^{\infty}_c$ be a cut-off of the ball of radius $\eta / 2$ around the origin. We write
\[\omega_0 = \omega_0\, \varphi + \omega_0\,(1- \varphi) \equiv \omega_0^F + \omega_0^N.\]      
We begin by observing that $\omega_0^N$ is a smooth function with compact support and hence 
\[\nabla \omega_{\nu} = \nabla H_{\nu} \ast \omega_0^F + H_{\nu} \ast \nabla \omega_0^N.\] Clearly  
$H_{\nu} \ast \nabla \omega_0^N$ is a bounded function, uniformly in $\nu$. Next we estimate  
$\nabla H_{\nu} \ast \omega_0^F$ far from the origin. Let $x$ be such that $|x|>\eta$. Then: 
\[ 
|\nabla H_{\nu} \ast \omega_0^F (x) | \leq \frac{1}{8\pi(\nu^2t^2)}\left| \int_{|y|<\eta/2} |x-y| e^{-|x-y|^2/(4\nu t)} \omega_0(y) \, dy \right| \] 
\[\leq \frac{C_1}{\nu^{3/2}t^{3/2}}\int_{|y|< \eta/2} e^{-C_2|x-y|^2/(\nu t)} \frac{1}{|y|} \, dy,\] 
where we have used the fact that there exist constants $C_1,\,C_2>0$ such that $|z|e^{-|z|^2} \leq C_1 e^{-C_2|z|^2}$, 
\[\leq \frac{2\pi C_1}{\nu^{3/2}t^{3/2}}e^{-C_2\eta/(\nu t)}\int_0^{\eta /2} e^{-r/(\nu t)} \, dr \leq C,\] 
for some $0<C<\infty$, $C$ independent of $\nu$. In summary we have shown that $\nabla \omega_{\nu}$ is  bounded 
in the complement of $B(0;\eta)$ uniformly in $\nu$.  
In view of this fact, since $Z^{\nu}(\omega_{\nu}) = \nu |\nabla \omega_{\nu}|^2$,  
\eqref{e:claim1} follows. This concludes the proof of the Claim. \end{proof}  
The proof of Theorem \ref{t:example} is complete.  
\end{proof}   

We close by verifying  that the sequence $\omega_{\nu}$ satisfies the 
hypothesis of Eyink's conjecture. 
 
\begin{proposition} \label{Besov} 
We have that $\omega_0 \in B^0_{2,\infty}$ is a viscosity solution of the 2D Euler equations and  
$\omega_{\nu} \rightharpoonup \omega_0$  weak-$\ast$ in
$L^{\infty}((0,T);B^0_{2,\infty})$, as $\nu \to 0^+$.
\end{proposition} 
 
\begin{proof} 
We begin by recalling the definition of the norm in  $B^0_{2,\infty}$ 
( see e.g. \cite{Trieb}, page 17):
\[
    \|f\|_{B^{0}_{2,\infty}} = \sup_{j\geq 0} \|\psi_{j}\ast f\|_{L^{2}},
\]
where $\psi_{j}$ are functions forming a Littlewood-Paley partion of
unity. In particular, the Fourier transform of $\psi_{0}$, 
$\widehat{\psi_{0}}$, is smooth, compactly supported in  the disk $B(0;1)$,
$\widehat{\psi_{0}}\equiv 1$ on  $B(0;2/3)$, while for $j>0$,
$\psi_{j}(x)= 2^{2j}\psi(2^{j} x)$, for a function $\psi$ such that
its Fourier transform $\widehat{\psi}$ is smooth, compactly supported
in the shell $\{1/2<|\xi|<2\}$, $\widehat{\psi}\equiv 1$ on  
$\{2/3<|\xi|<4/3\}$. 

We will estimate the low and high-frequency contribution to the
$B^{0}_{2,\infty}$-norm of $\omega_{\nu}= H_{\nu}\ast \omega_{0}$
separately. Here again, $H_{\nu}$ is the heat kernel and the
convolution is only in the space variable.

For the low-frequency part, we  observe that
$\omega_{0} \in L^{1}_{c}(\real^{2})$, the \mbox{$L^{1}$-norm} of
$H_{\nu}$ as a function of $x$ is uniformly bounded in $t$ and $\nu$,
and that $\psi_{0}$ is smooth, rapidly decreasing. Consequently, by Young's
inequality
\begin{equation} \label{e:norm1}
   \begin{aligned}
   \|\psi_{0}\ast \omega_{\nu}(t)\|_{L^{2}} &\leq
    \|\psi_{0}\|_{L^{2}}\,\|\omega_{\nu}(t)\|_{L^{1}} \\
     &\leq  \|\psi_{0}\|_{L^{2}}\,
     \|H_{\nu}(t)\|_{L^{1}}\,\|\omega_{0}\|_{L^{1}} \leq C,
 \end{aligned}
\end{equation}
$C$ independent of $\nu$ and $t$.

To bound the high-frequency part we will
employ the Fourier transform and knowledge of the behavior of
$\widehat{w_{0}}$ gained in Proposition \ref{e:normest}.
In view of \eqref{e:square}, we can write
\[
  \begin{aligned}
    \|\psi_{j}\ast \omega_{\nu}(t)\|^{2}_{L^{2}} &=
     \int_{\real^{2}} |\widehat{\psi_{j}}(\xi)|^2\,
     |\widehat{\omega_{\nu}}|^{2}\,d\xi = \int_{\real^{2}} 
     |\widehat{\psi}(2^{-j}\xi)|^2\,e^{-2 \nu t\,|\xi|^{2}}
     |\widehat{\omega_{0}}(\xi)|^{2}\,d\xi  \\
     &\leq \int_{\real^{2}} 
     |\widehat{\psi}(2^{-j}\xi)|^2\,
     |\widehat{\omega_{0}}(\xi)|^{2}\,d\xi  \\
     &= \int_{\real^{2}} 
     |\widehat{\psi}(2^{-j}\xi)|^2\,\frac{1}{|\xi|^{2}} (
     4\pi^2 + 4\pi e(\xi) + |e(\xi)|^2)\,d\xi .
 \end{aligned}
\] 
We now change variables from $\xi$ to 
$\xi'=2^{-j}\xi$, and use the support properties of $\widehat{\psi}$ to obtain
\begin{equation} \label{e:norm2}
     \|\psi_{j}\ast \omega_{\nu}\|^{2}_{L^{2}}\leq \int_{1/2<|\xi|<2} 
     |\widehat{\psi}(\xi')|^2\,\frac{1}{|\xi'|^{2}} (
     4\pi^2 + 4\pi e(2^{j} \xi') + |e(2^{j} \xi')|^2)\,d\xi' \leq C,
\end{equation}
with $C$ again independent of $\nu$ and $t$,
since the function $e(2^{j}\xi)$ is bounded uniformly in $j$ and $\xi$.
We remark that this also shows that $\omega_{0}\in
B^{0}_{2,\infty}$.

Combining \eqref{e:norm1} and  \eqref{e:norm2}  finally gives 
\[
      \sup_{\nu >0}
      \|\omega_{\nu}\|_{L^{\infty}((0,T);B^{0}_{2,\infty})} \leq C <\infty.
\]
Therefore, there exists a subsequence $\omega_{\nu_{k}}$, which converges  
weak-$\ast$ in  $L^{\infty}((0,T);B^{0}_{2,\infty})$, to a weak-$\ast$ limit. 
But, since $H_{\nu}\rightharpoonup \delta_{0}$ in $\Sc'$ as $\nu \to 0^+$, we conclude that the
whole family $\omega_{\nu}$ converges  weak-$\ast$ in
$L^{\infty}((0,T);B^{0}_{2,\infty})$ and the weak-$\ast$ limit is $\omega_0$.
\end{proof}

What we have actually accomplished with Theorem \ref{t:example} and Proposition \ref{Besov} is to give a counterexample to the part of Eyink's conjecture identifying viscous and transport enstrophy defects. We have answered in the affirmative the part of the conjecture regarding the existence of a nontrivial enstrophy defect. Although we found such an example only for the 
viscous enstrophy defect, this is the physically meaningful
one. Clearly, from the point of view of turbulence theory, one should
attempt to understand better the viscous enstrophy defect. Informally, viscous dissipation of a quantity is enhanced the more complicated the spatial distribution of that quantity. Our radially 
symmetric, monotonic example is as simple a configuration as possible, and, as such, should have 
the least dissipation. We imagine that, in some sense, the viscous enstrophy defect should be greater for a generic 
configuration, and existence of the viscous enstrophy defect would be the more problematic issue.  

\section{Conclusions} 
 
        We would like to add a few general remarks regarding the work presented here. First, the theory 
of viscous and transport enstrophy defects can be formulated in the more general setting of weak and renormalized solutions of linear transport equations and vanishing viscosity limits. The only instance where the specific form of the  incompressible fluid flow equations was used is when we attributed meaning to the expression  
$u|\omega|^2$ for $\omega \in L^2(\log L)^{1/4}$. In particular, the counterexample presented in Section 6 is 
really a solution of the heat equation, of some interest even without mentioning the fluid dynamical context.     
 
        Our counterexample to Eyink's Conjecture is circularly symmetric, and as such, it corresponds to solutions of the Navier-Stokes equations for which the nonlinear term $P(u \cdot \nabla u)$ vanishes identically ($P$ is the Leray projector). Since turbulence is regarded as coming from the interaction of nonlinearity and small viscosity, it is fair to ask what possible relevance would such an example have for the understanding of turbulence. If one looks at the cascade ansatz, the basic idea is that the nonlinearity produces a flow of enstrophy, from large to small scales across the inertial range, to be dissipated by viscosity. For flows with finite enstrophy, the nonlinearity must play a crucial role in sustaining the cascade because without the nonlinearity the viscosity would instantly make small scale enstrophy disappear. Now, for flows with infinite enstrophy, the nonlinearity is not needed for a sustained cascade because there already is an infinite supply of enstrophy at small scales. At this level, it is possible for the flow of enstrophy to small scales due to the nonlinearity to be small, or irrelevant. This would be a plausible explanation for 
why the viscous and transport enstrophy defects are not the same. It would be interesting to take a new look at the 
Kraichnan-Batchelor theory in light of this possibility.    
 
        It is not clear whether the notions of enstrophy defect will become useful in general issues of interest in PDE, but this is certainly possible and further research along this line is amply warranted. 
Due to the unexplored nature of this subject, it is easy to formulate a long list of  open problems. 
We will single out a few that appear either particularly accessible or interesting. The main open problem is to prove that viscous enstrophy defects are well defined for some class of flows with infinite (local) enstrophy. Another important problem is to find an example of a solution to an inviscid transport equation, preferably given by a solution of the Euler equations, for which the transport enstrophy defect is nonzero. We have seen that the transport enstrophy density is a weak solution of the appropriate transport equation for initial vorticities in  
$L^2(\log L)^{1/4}$ if the weak solution comes from vanishing viscosity. It would be very interesting to find other properties of viscosity solutions that are not shared by general weak solutions. Although enstrophy plays a distinguished role among integrals of convex functions of vorticity due to its relevance to turbulence modeling, it is  reasonable to ask to which extent similar defects might be usefully associated to other such first integrals. There is a certain arbitrariness in the definition of transport enstrophy defect that might be explored, as one could define another inviscid enstrophy defect by mollifying the initial data, for example. Finally, we state again a problem suggested in Section 2: determine whether viscosity solutions are renormalized solutions of the transport equations if initial vorticity is in $L^p$, $p<2$. Note that nonuniqueness of weak solutions follows immediately if this is not the case.  
   
{\small Acknowledgments: Research of M.C. Lopes Filho is supported in
  part by CNPq grant \#300.962/91-6. Research of H. J. Nussenzveig Lopes is
  supported in part by CNPq grants \#300.158/93-9 and \#200.951/03-3 . This
  research has been supported by FAPESP grant \#02/05556-2. This work was  
  partially conducted while M.C. Lopes Filho and H.J. Nussenzveig Lopes were on sabbatical 
  leave at Penn State University. The third author would
  like to thank C. Doering and E. Titi, for their comments
  and insight. The second author would like to thank R. Beals, R. Coifman, and  P. Jones for helpful discussions. She
  would like to acknowledge also the continuous interest of her
  Ph.D. adviser M. E. Taylor. This work was partially conducted while
  A. L. Mazzucato was a Gibbs Instructor at Yale University.}


\begin{thebibliography}{LNT00}

\bibitem[Ada75]{Adams}
R.~A. Adams, \emph{Sobolev spaces}, Academic Press [A subsidiary of Harcourt
  Brace Jovanovich, Publishers], New York-London, 1975, Pure and Applied
  Mathematics, Vol. 65.

\bibitem[Bat69]{Ba}
 G.~K. Batchelor, \emph{Computation of the energy spectrum in homogeneous
  two-dimensional turbulence}, Phys. Fluids Suppl. II \textbf{12} (1969),
  233--239.


\bibitem[BL76]{B-L}
J.~Bergh and J.~L{\"o}fstr{\"o}m, \emph{Interpolation spaces. {A}n
  introduction}, Springer-Verlag, Berlin, 1976, Grundlehren der Mathematischen
  Wissenschaften, No. 223.

\bibitem[BS88]{BS88}
C. Bennett and R. Sharpley, \emph{Interpolation of operators}, Pure and
  Applied Mathematics, vol. 129, Academic Press Inc., Boston, MA, 1988.

\bibitem[Del91]{DeL}
J.-M. Delort, \emph{Existence de nappes de tourbillon en dimension deux}, J.
  Amer. Math. Soc. \textbf{4} (1991), no.~3, 553--586.

\bibitem[DL89]{DiPL}
R.~J. DiPerna and P.-L. Lions, \emph{Ordinary differential equations, transport
  theory and {S}obolev spaces}, Invent. Math. \textbf{98} (1989), no.~3,
  511--547.

\bibitem[DM87]{DiPM1}
R.~J. DiPerna and A.~J. Majda, \emph{Concentrations in regularizations for
  {$2$}-{D} incompressible flow}, Comm. Pure Appl. Math. \textbf{40} (1987),
  no.~3, 301--345.

\bibitem[DR00]{DR}
J.~Duchon and R.~Robert, \emph{Inertial energy dissipation for weak solutions
  of incompressible {E}uler and {N}avier-{S}tokes equations}, Nonlinearity
  \textbf{13} (2000), no.~1, 249--255.

\bibitem[Eva90]{EvansCBMS}
L.~C. Evans, \emph{Weak convergence methods for nonlinear partial differential
  equations}, CBMS Regional Conference Series in Mathematics, vol.~74,
  Published for the Conference Board of the Mathematical Sciences, Washington,
  DC, 1990.

\bibitem[Eyi01]{Ey}
G.~L. Eyink, \emph{Dissipation in turbulent solutions of 2{D} {E}uler
  equations}, Nonlinearity \textbf{14} (2001), no.~4, 787--802.

\bibitem[Fri95]{Frisch}
U.~Frisch, \emph{Turbulence}, Cambridge University Press, Cambridge, 1995, The
  legacy of A. N. Kolmogorov.

\bibitem[Kra67]{Kraich}
R.~H. Kraichnan, \emph{Inertial ranges in two-dimensional turbulence}, Phys
  Fluids \textbf{10} (1967), 1417--1423.


\bibitem[Lio96]{Lio}
P.-L. Lions, \emph{Mathematical topics in fluid mechanics. {V}ol. 1}, Oxford
  Lecture Series in Mathematics and its Applications, vol.~3, The Clarendon
  Press Oxford University Press, New York, 1996, Incompressible models, Oxford
  Science Publications.

\bibitem[LNT00]{LLT}
M.~C. Lopes, Filho, H.~J. Nussenzveig, Lopes, and E.~Tadmor, \emph{Approximate
  solutions of the incompressible {E}uler equations with no concentrations},
  Ann. Inst. H. Poincar\'e Anal. Non Lin\'eaire \textbf{17} (2000), no.~3,
  371--412.

\bibitem[Maj93]{Maj}
A.~J. Majda, \emph{Remarks on weak solutions for vortex sheets with a
  distinguished sign}, Indiana Univ. Math. J. \textbf{42} (1993), no.~3,
  921--939.

\bibitem[MB02]{majbert}
A.~J. Majda and A.~L. Bertozzi, \emph{Vorticity and incompressible flow},
  Cambridge Texts in Applied Mathematics, vol.~27, Cambridge University Press,
  Cambridge, 2002.

\bibitem[Sch95]{Sch}
S. Schochet, \emph{The weak vorticity formulation of the {$2$}-{D} {E}uler
  equations and concentration-cancellation}, Comm. Partial Differential
  Equations \textbf{20} (1995), no.~5-6, 1077--1104.

\bibitem[Shn00]{Shni}
A.~Shnirelman, \emph{Weak solutions with decreasing energy of incompressible
  {E}uler equations}, Comm. Math. Phys. \textbf{210} (2000), no.~3, 541--603.

\bibitem[Ste70]{St70}
E.~M. Stein, \emph{Singular integrals and differentiability properties of
  functions}, Princeton Mathematical Series, No. 30, Princeton University
  Press, Princeton, N.J., 1970.

\bibitem[Tri92]{Trieb}
H. Triebel, \emph{Theory of function spaces. {II}}, Monographs in
  Mathematics, vol.~84, Birkh\"auser Verlag, Basel, 1992.

\bibitem[Vis98]{Vi1}
M.~Vishik, \emph{Hydrodynamics in {B}esov spaces}, Arch. Ration. Mech. Anal.
  \textbf{145} (1998), no.~3, 197--214.

\bibitem[Vis99]{Vi2}
M.~Vishik, \emph{Incompressible flows of an ideal fluid with vorticity in
  borderline spaces of {B}esov type}, Ann. Sci. \'Ecole Norm. Sup. (4)
  \textbf{32} (1999), no.~6, 769--812.

\bibitem[VW93]{VW93}
I.~Vecchi and S.~J. Wu, \emph{On {$L\sp 1$}-vorticity for {$2$}-{D}
  incompressible flow}, Manuscripta Math. \textbf{78} (1993), no.~4, 403--412.

\bibitem[Yud63]{Yu63}
V.~I. Yudovi{\v{c}}, \emph{Non-stationary flows of an ideal incompressible
  fluid}, \u Z. Vy\v cisl. Mat. i Mat. Fiz. \textbf{3} (1963), 1032--1066.

\bibitem[Yud95]{Yu}
V.~I. Yudovich, \emph{Uniqueness theorem for the basic nonstationary problem in
  the dynamics of an ideal incompressible fluid}, Math. Res. Lett. \textbf{2}
  (1995), no.~1, 27--38.

\end{thebibliography}
 
\def\ocirc#1{\ifmmode\setbox0=\hbox{$#1$}\dimen0=\ht0 \advance\dimen0
  by1pt\rlap{\hbox to\wd0{\hss\raise\dimen0
  \hbox{\hskip.2em$\scriptscriptstyle\circ$}\hss}}#1\else {\accent"17 #1}\fi}
  \def\cprime{$'$}
\providecommand{\bysame}{\leavevmode\hbox to3em{\hrulefill}\thinspace}

\end{document}